\documentclass[a4paper,10pt,leqno]{amsart}
\title{Twisted $L^2$-torsion functions  obstruct  fibring over the circle}
\author{Sam Hughes}
\address{Max-Planck-Institut für Mathematik\\ 
Vivatsgasse 7\\ 
53111 Bonn\\ 
Germany}
              \email{sam.hughes.maths@gmail.com}\email{hughess@mpim-bonn.mpg.de}
              \urladdr{https://samhughesmaths.github.io}
\author{Wolfgang L\"uck}
        \address{Mathematisches Institut der Universit\"at Bonn\\
                Endenicher Allee 60\\
                53115 Bonn, Germany}
         \email{wolfgang.lueck@him.uni-bonn.de}\email{lueck@math.uni-bonn.de}
         \urladdr{http://www.him.uni-bonn.de/lueck}
         \date{November 2026}
         \keywords{fibring over the circle, $L^2$-torsion functions and their asymptotic behaviour}
         \makeatletter
     \@namedef{subjclassname@2020}{%
  \textup{2020} Mathematics Subject Classification}
\makeatother
\subjclass[2020]{55R10 (primary); 57Q10, 58J52, 46L99 (secondary)}

 \usepackage[utf8]{inputenc}
    \usepackage{comment}
  \usepackage{calc}
  \usepackage{enumerate,amssymb,bm}
  \usepackage[arrow,curve,matrix,tips,2cell]{xy}
    \SelectTips{eu}{10} \UseTips
    \UseAllTwocells
  \usepackage{tikz}
  \usetikzlibrary{decorations.pathreplacing}
  \usepackage{color}
    \usepackage{scalerel}
  \usepackage{braket}
   \usepackage{mathtools}
  \usepackage{bbm}
  \usepackage{enumitem}
  \usepackage{float}
   \DeclareMathAlphabet{\matheurm}{U}{eur}{m}{n}
\usepackage[pagebackref=true]{hyperref}
  \usepackage{datetime2}

\renewcommand*{\backref}[1]{}
\renewcommand*{\backrefalt}[4]
{\ifcase #1
        No citation in the text.
    \or
        Cited on Page #2.
    \else
        Cited on Pages #2.
    \fi
}

\makeatletter
\def\@tocline#1#2#3#4#5#6#7{\relax
  \ifnum #1>\c@tocdepth 
  \else
    \par \addpenalty\@secpenalty\addvspace{#2}%
    \begingroup \hyphenpenalty\@M
    \@ifempty{#4}{%
      \@tempdima\csname r@tocindent\number#1\endcsname\relax
    }{%
      \@tempdima#4\relax
    }%
    \parindent\z@ \leftskip#3\relax \advance\leftskip\@tempdima\relax
    \rightskip\@pnumwidth plus4em \parfillskip-\@pnumwidth
    #5\leavevmode\hskip-\@tempdima
      \ifcase #1
       \or\or \hskip 1em \or \hskip 2em \else \hskip 3em \fi%
      #6\nobreak\relax
    \hfill\hbox to\@pnumwidth{\@tocpagenum{#7}}\par
    \nobreak
    \endgroup
  \fi}
\makeatother

\DeclareMathAlphabet{\matheurm}{U}{eur}{m}{n}

\DeclareMathOperator{\cone}{cone}

\DeclareMathOperator{\el}{\operatorname{el}}

\DeclareMathOperator{\Ext}{Ext}
\DeclareMathOperator{\F}{F}

\DeclareMathOperator{\FP}{FP}

\DeclareMathOperator{\GL}{GL}

\DeclareMathOperator{\id}{id}

\DeclareMathOperator{\lead}{lead}

\DeclareMathOperator{\ord}{ord}

\DeclareMathOperator{\pr}{pr}

\DeclareMathOperator{\tors}{tors}
\DeclareMathOperator{\tr}{tr}

\DeclareMathOperator{\Wh}{Wh}
\DeclareMathOperator{\width}{width}
\DeclareMathOperator{\zero}{zero}

\newcommand{\UFD}{\textup{UFD}}

  \newcommand{\IC}{\mathbb{C}}

  \newcommand{\IF}{\mathbb{F}}

  \newcommand{\IQ}{\mathbb{Q}}
  \newcommand{\IR}{\mathbb{R}}

  \newcommand{\IZ}{\mathbb{Z}}

  \newcommand{\calb}{\mathcal{B}}

  \newcommand{\cald}{\mathcal{D}}
  \newcommand{\cale}{\mathcal{E}}

  \newcommand{\caln}{\mathcal{N}}

\newcommand{\GLpar}[2]{\GL_{#1}(#2)}

\newcommand{\Mpar}[3]{\operatorname{M}_{#1,#2}(#3)}

\newcounter{commentcounter}

\theoremstyle{plain}
\newtheorem{theorem}{Theorem}[section]

\newtheorem{lemma}[theorem]{Lemma}
\newtheorem{corollary}[theorem]{Corollary}

\newtheorem*{theorem*}{Theorem}
\newtheorem*{theoremA*}{Theorem A}
\newtheorem*{theoremB*}{Theorem B}

\theoremstyle{definition}
\newtheorem{definition}[theorem]{Definition}

\newtheorem{example}[theorem]{Example}

\newtheorem{remark}[theorem]{Remark}

\newtheorem{basicsetup}[theorem]{Basic Setup}

\newtheorem*{definition*}{Definition}

\theoremstyle{remark}

\makeatletter\let\c@equation=\c@theorem\makeatother

\theoremstyle{definition}

\newcounter{othercommentcounter}

 \newcommand{\version}[1] 
     {\begin{center} last edited on #1\\
         last compiled on \today \ at \DTMcurrenttime.\\
         name of tex-file: \jobname
       \end{center}}

\newcommand{\RALI}{\textup{RALI}}
     \newcommand{\RFRS}{\textup{RFRS}}

\begin{document}

     \begin{abstract}
       The vanishing of the $L^2$-Betti numbers of the universal covering of a closed
       manifold is an obstruction to fibring over the circle.  When the $L^2$-Betti numbers vanish
       we introduce a secondary obstruction based on the asymptotic behaviour of the
       twisted $L^2$-torsion function at $0$ and $\infty$.  We provide a number of examples 
       showing the invariant is indeed non-trivial.
       \end{abstract}

\maketitle

\newlength{\origlabelwidth} \setlength\origlabelwidth\labelwidth


\typeout{------------------- Section 1: Introduction -----------------}

\section{Introduction}\label{sec:introduction}

A compact manifold $M$ \emph{fibres over the circle} if we may view it as a fibre bundle
$F\to M\to S^1$ where $F$ is a compact manifold.  Two major results concern hyperbolic
manifolds.  Firstly Agol's positive resolution~\cite{Agol(2013)} to Thurston's
question~\cite[Page~380]{Thurston(1982)} whether every closed hyperbolic $3$-manifold
virtually fibres over the circle $S^1$.  Secondly Italiano--Martelli--Migliorini fibred a
cusped hyperbolic $5$-manifold over the circle~\cite{Italiano-Martelli-Migliorini(2023five)}.

There are
few obstructions to fibring a closed manifold over the circle.  One of the main
ones being L\"uck's theorem~\cite{Lueck(1993)} on non-vanishing $L^2$-Betti numbers, as
explained in Section~\ref{sec:L2-Betti_numbers_are_obstruction_to_fibring_over_S1}.
Kielak~\cite{Kielak(2020fibring)} proved that for finitely generated \RFRS-groups the
vanishing of the first $L^2$-Betti number is the only obstruction to the group virtually
admitting a homomorphism to $\IZ$ with finitely generated kernel, see
Theorem~\ref{the:Kielak}. Fisher~\cite{Fisher(2024)} extended Kielak's result to higher
dimensions.  Namely a \RFRS\ group of type $\FP_n(\IQ)$ virtually admits a homomorphism to
$\IZ$ with kernel of type $\FP_n(\IQ)$ if and only if all the $L^2$-Betti numbers in
dimensions $\le n$ vanish, see Theorem~\ref{the:Fisher_Q}.  
In~\cite{Hughes-Leary-Lueck(2026not_fibre)} examples of closed manifolds with \RFRS\
fundamental group and vanishing $L^2$-Betti numbers that do not virtually fibre over the
circle were constructed.

We will present a new secondary obstruction for fibring in
Theorem~\ref{the:Property_special_is_an_obstruction_to_fibring_and_virtually_fibring}.
It deals with the \emph{reduced $L^2$-torsion function}
\[
  \rho^{(2)}(\widetilde{X};\phi) \in \cale
\]
of an $L^2$-acyclic connected finite $CW$-complex $X$, see Section~\ref{sec:prelims} for the definitions.  
It is secondary in the sense that
it is only defined if the first condition, the vanishing of all $L^2$-Betti numbers, is
satisfied.


\begin{theorem}[Theorem~\ref{the:Property_special_is_an_obstruction_to_fibring_and_virtually_fibring}]%
  \label{the:Property_special_is_an_obstruction_to_fibring_and_virtually_fibring_intro}
  Let $X$ be a connected finite $CW$-complex satisfying one of the following equivalent
  conditions:

  \begin{enumerate}

  \item There is a connected finite $CW$-complex $F$ such that $X$ is homotopy equivalent
    to the mapping torus $T_f$ for some selfhomotopy equivalence
    $f \colon F \xrightarrow{\simeq} F$;

  \item There exists a fibration $F \xrightarrow{i} E \xrightarrow{p} S^1$ such that $E$
    is homotopy equivalent to $X$ and $F$ is a connected finite $CW$-complex.

  \end{enumerate}

  Then:
  \begin{enumerate}
  
  \item\label{the:mapping_tori:L2-acyclic_intro} The $L^2$-Betti number
    $b_n^{(2)}(\widetilde{X})$ vanishes for every $n \in \IZ_{\ge 0}$;

  \item\label{the:mapping_tori:overline(rho)_intro} Let $\phi \colon \pi_1(X) \to \IZ$ be the
    surjective group homomorphism which is either associated to $T_f$
    or comes from $\pi_1(p)$ for the projection $p \colon E \to S^1$. Suppose that $\pi_1(X)$
    is residually finite and satisfies the Farrell--Jones Conjecture.

    Then there are elements $T_0, T_{\infty} \in \IR_{> 0}$  and a
    representative $\rho \colon \IR_{> 0} \to \IR$ of
    $\overline{\rho}(\widetilde{X};\phi)$ satisfying
    \[
      \rho(t) = \begin{cases} 0 & \text{for} \; t < T_0;
        \\
        \chi(F) \cdot \ln(t) & \text{for} \; T_{\infty} < t.
      \end{cases}
    \]
    In particular $\overline{\rho}(\widetilde{X};\phi)$ is special.
  \end{enumerate}
\end{theorem}

Here the $\phi$-twisted $L^2$-torsion is
\emph{special} if there is a representative
satisfying the conclusion of the theorem, that is there are elements
$T_0, T_{\infty} \in \IR_{> 0}$, an integer $k$,  and a representative
$\rho \colon \IR_{> 0} \to \IR$ of $\overline{\rho}(\widetilde{X};\phi)$ satisfying
    \[
      \rho(t) = \begin{cases} 0 & \text{for} \; t < T_0;
        \\
        k\cdot \ln(t) & \text{for} \; T_{\infty} < t.
      \end{cases}
    \]

    We show that the obstruction is non-trivial by computing the $\phi$-twisted $L^2$-torsion
    function for integral characters of right-angled Artin groups.  Taking $L$ to be a
    minimal flag triangulation of the real projective plane 
    gives examples of characters
    that are not special, see Example~\ref{ex:RAAG_RP2}.

    \begin{theorem}[Theorem~\ref{thm:phi_twisted_torsion_RAAGs}]%
\label{thm:phi_twisted_torsion_RAAGs_intro}
Let $L$ be a $\IQ$-acyclic flag complex, let $A_L$ be the right-angled Artin group on $L$,
and let $\phi:A_L \to \IZ$ be a group homomorphism.  Define two constants
\begin{eqnarray*}
  C_L
  & = &
\sum_{i\geq 2}(-1)^i \cdot \ln(|\mathrm{tor} H_{i-1}(L)|);
\\
  C_\phi
  & =&
\sum_{v\in L^{(0)}} |\phi(v)| \cdot \overline \chi(\mathrm{Lk}_{L}(v)).
\end{eqnarray*}
Then a representative $\rho$ for $\overline\rho^{(2)}(A_L;\phi)(t)$ is given by 
\[
t\mapsto
\begin{cases} 
  C_L- C_\phi \cdot \ln(t) & t\geq 1;
  \\
C_L & t \le 1. 
\end{cases}
\]
\end{theorem}

We consider a number of other examples.  In
Section~\ref{sec:The_special_case_of_an_infinite_cyclic_group} we
consider the special case of $G=\IZ$ and $\phi = \id_{\IZ}$ in detail
and we relate the $\phi$-twisted torsion function to Alexander
polynomials.   In
Subsection~\ref{subsec:Aspherical} we
revisit an example from~\cite{Hughes-Leary-Lueck(2026not_fibre)} and
show that its $\phi$-twisted $L^2$-torsion function is special. 
  In Subsection~\ref{subsec:poly_growth} we
compute the $\phi$-twisted $L^2$-torsion function of a semidirect
product $\Gamma=G\rtimes_\Theta \IZ$ where $\Gamma$ is torsion-free
residually finite, $G$ is a hyperbolic relative to a collection of
virtually polycyclic groups, and $\Theta$ is polynomially growing.
The reader should compare this result with~\cite[Theorem~7.19]{Hughes-Lueck(2025)}.
 Finally, in Subsection~\ref{subsec:The_matrix_3g_plus_2e} we compute
 the $\phi$-twisted $L^2$-torsion function
of the operator associated to the matrix $3g+2e$ and show it is not
special.  This is particularly relevant to the counterexamples
constructed in~\cite{Hughes-Leary-Lueck(2026not_fibre)}.


\subsection*{Structure of the paper}
The paper is organized as follows:
\tableofcontents


\subsection*{Acknowledgments}\label{subsec:Acknowledgements}
SH was supported by a Humboldt Research Fellowship at Universit\"at Bonn 
and by the Max Planck Institute for Mathematics in Bonn.
The paper is funded by the Deutsche Forschungsgemeinschaft (DFG, German Research
Foundation) under Germany's Excellence Strategy \--- GZ 2047/1, Projekt-ID 390685813,
Hausdorff Center for Mathematics at Bonn.


\subsection*{AI declaration} The authors did not use AI in the
preparation of this article.


\typeout{----------------------- Section 2: Preliminaries ---------------------------}

\section{Preliminaries}\label{sec:prelims}


\subsection{Basics about mapping tori of chain complexes}%
\label{sec:Basics_about_mapping_tori_of_chain_complexes}

For a connected finite $CW$-complex $X$ and a group epimorphism
$\phi \colon \pi_1(X) \to \IZ$ the following assertions are known to be equivalent:

\begin{enumerate}
\item There is a fibration $F \to X \xrightarrow{p} S^1$ such that $F$ is a connected
  finite $CW$-complex and $\pi_1(p) = \phi$;

\item There is a connected finite $CW$-complex $F$, a selfhomotopy equivalence
  $f \colon F \xrightarrow{\simeq} F$, and a homotopy equivalence $u \colon T_f \to X$
  with the mapping torus of $f$ as source such that
  $ \phi \circ \pi_1(u) \colon \pi_1(T_f) \to \pi_1(S^1) = \IZ$ is $\pi_1(\pr)$ for the
  canonical projection $\pr \colon T_f \to S^1$;

\item The total space $\overline{X}$ of the infinite cyclic covering
  $q \colon \overline{X} \to X$ associated to $\phi$ is homotopy equivalent to a connected
  finite $CW$-complex.

\end{enumerate}

Next we explain the chain complex version of this fact.

Let $\Psi \colon S \xrightarrow{\cong} S$ be a ring automorphism.  Denote by
$S_{\Psi}[t,t^{-1}] = S_{\Psi}[\IZ]$ the $\Psi$-twisted group ring of $\IZ$ with
coefficients in $S$.  Let $D_*$ be an $S$-chain complex.  Let
$f_* \colon D_* \to \Psi^*D_*$ be an $S$-chain map, where $\Psi^*D_*$ is the $S$-chain
complex obtained from $D_*$ by restriction with $\Psi$. We obtain an
$S_{\Psi}[t,t^{-1}]$-chain map
\[
  \widehat{f_*} \colon S_{\Psi}[t,t^{-1}] \otimes_{S} D_* \to S_{\Psi}[t,t^{-1}]
  \otimes_{S} D_*
\]
by sending $r \otimes x$ for $r \in S_{\Psi}[t,t^{-1}]$ and $x \in D_n$ to
$r \otimes x - rt^{-1} \otimes f_n(x)$.

In the sequel we call a chain complex \emph{finite} if it is finitely generated and finite-dimensional.

\begin{definition}[Mapping torus for chain
  maps]\label{def:chain_complex_version_of_a_mapping_torus}
  Define the mapping torus $T(f_*)_*$ of $f_*$ to be the $S_{\Psi}[t,t^{-1}]$-chain
  complex given by the mapping cone $\cone(\widehat{f_*})$ of $\widehat{f}_*$.
\end{definition}

\begin{lemma}\label{lem:elementary_property_of_mapping_torus}
  The following assertions are equivalent for a finite free $S_{\Psi}[t,t^{-1}]$-chain
  complex $C_*$:

  \begin{enumerate}
  \item\label{lem:elementary_property_of_mapping_torus:(1)} The $S$-chain complex $i^*C_*$
    is $S$-chain homotopy equivalent to a finite free $S$-chain complex, where $i^*C_*$ is
    obtained from $C_*$ by restriction with the inclusion
    $i \colon S \to S_{\Psi}[t,t^{-1}]$;

  \item\label{lem:elementary_property_of_mapping_torus:(2)} There exists a finite free
    $S$-chain complex $D_*$ and an $S$-chain homotopy equivalence
    $f_* \colon D_* \to \Psi^*D_*$ such that $C_*$ is $S_{\Psi}[t,t^{-1}]$-chain homotopy
    equivalent to $T(f_*)_*$.
  \end{enumerate}

  If one replaces ``finite free'' by ``finite projective'', ``finitely generated free'',
  or ``finitely generated projective'' above, the equivalence above remains  true.

\end{lemma}
\begin{proof} This follows from~\cite[Lemma~3.3]{Hughes-Leary-Lueck(2026not_fibre)}.
\end{proof}


\subsection{\texorpdfstring{$L^2$}{L2}-Betti numbers are obstruction to fibring over \texorpdfstring{$S^1$}{S1}}%
\label{sec:L2-Betti_numbers_are_obstruction_to_fibring_over_S1}

$L^2$-Betti numbers are obstructions to fibring over $S^1$ as explained next.  Let
$F \to E \to S^1$ be a fibration such that $F$ has the homotopy type of a $CW$-complex of
finite type.  Then $E$ has the homotopy type of a $CW$-complex of finite type and is
actually homotopy equivalent to a mapping torus $T_f$ of a selfmap $f \colon Y \to Y$ of
some $CW$-complex $Y$ which is of finite type and homotopy equivalent to $F$.  This does
imply that the $L^2$-Betti numbers $b_n^{(2)}(\widetilde{E})$ vanish for $n \in \IZ_{\ge 0}$,
see~\cite[Theorem~6.63 on page~270]{Lueck(2002)}.

An analogous statement holds also for the $\IF_p$-version of the $L^2$-Betti numbers
$b_n(\widetilde{E};\cald_{\IF_p[\pi_1(E)]})$,
see~\cite[Theorem~3.25]{Avramidi-Lueck(2026)}. These numbers are defined based on the work
of Jaikin-Zapirain and systematically studied by Avramidi--L\"uck~\cite{Avramidi-Lueck(2026)},
provided that $\pi_1(E)$ is a \RALI-group, i.e., is
residually (amenable and locally indicable).

Note that the class of finitely generated RALI-groups is larger than the class of finitely
generated RFRS-groups, since a finitely generated RFRS-group is residually (locally
indicable and virtually abelian), see~\cite[Section~6]{Okun-Schreve(2024orders)}.

In the sequel $\caln(G)$ denotes the group von Neumann algebra which can be identified with
the algebra $\calb(L^2(G),L^2(G))^G$ of bounded $G$-equivariant operators $L^2(G) \to L^2(G)$.

A group $G$ is \emph{algebraically fibreed} if it admits a homomorphism
$\phi \colon G \to \IZ$ whose kernel is finitely generated.  Next we record the celebrated
result of Kielak~\cite[Theorem~5.3]{Kielak(2020fibring)}.

\begin{theorem}\label{the:Kielak}
  Let $G$ be an infinite finitely generated group which is virtually an $\RFRS$-group.

  Then $G$ virtually algebraically fibres if and only if its first $L^2$-Betti number
  $b_1^{(2)}(G) = b_1^{(2)}(EG;\caln(G))$ vanishes.
\end{theorem}

This was extended by Fisher~\cite[Theorem~A]{Fisher(2024)} to higher dimensions.

\begin{theorem}\label{the:Fisher_Q}
  Let $G$ be a group which is virtually an $\RFRS$-group. Consider $n \in \IZ_{\ge 1}$.
  Suppose that $G$ is of type $\FP_n(\IQ)$. Then the following statements are equivalent:

  \begin{enumerate}
  \item\label{the:Fisher:(1)} We have $b_p^{(2)}(G) = 0$ for $p \le n$;
  \item\label{the:Fisher_Q:(2)} There is a finite-index subgroup $H\subseteq G$ and an
    epimorphism $\phi \colon H \to \IZ$ such that $\ker(\phi)$ is of type $\FP_n(\IQ)$;

  \item\label{the:Fisher_Q:(3)} There is a finite-index subgroup $H' \subseteq G$ and an
    epimorphism $\phi' \colon H' \to \IZ$ such that $b_p(H';\IQ) < \infty $ holds for
    $p \le n$.
  \end{enumerate}
\end{theorem}

There is also the following version for an arbitrary field $F$ due to
Fisher~\cite[Theorem~B]{Fisher(2024)}.

 \begin{theorem}\label{the:Fisher_F}
   Let $G$ be a group which is virtually an $\RFRS$-group. Consider $n \in \IZ_{\ge 1}$.
   Let $F$ be a field. Suppose that $G$ is of type $\FP_n(F)$. Then the following
   statements are equivalent:

   \begin{enumerate}
   \item\label{the:Fisher_F:(1)} We have $b_p(G;\cald_{F[G]}) = 0$ for $p \le n$;
   \item\label{the:Fisher_F:(2)} There is a finite-index subgroup $H\subseteq G$ and an
     epimorphism $\phi \colon H \to \IZ$ such that $\ker(\phi)$ is of type $\FP_n(F)$;

   \item\label{the:Fisher_F:(3)} There is a finite-index subgroup $H' \subseteq G$ and an
     epimorphism $\phi' \colon H' \to \IZ$ such that $b_p(H';F)< \infty$ holds for
     $p \le n$.
   \end{enumerate}
 \end{theorem}

 Note that one cannot replace in Theorems~\ref{the:Kielak},~\ref{the:Fisher_Q},
 and~\ref{the:Fisher_F} the condition $\RFRS$-group by the weaker condition that the group
 is residually (torsionfree and nilpotent),
 see~\cite[Remark~5.8]{Hughes-Leary-Lueck(2026not_fibre)}.


 \subsection{{Review of \texorpdfstring{$\phi$}{phi}-twisted \texorpdfstring{$L^2$}{L2}-torsion function for chain complexes}}%
 \label{sec:Review_of_varphi-twisted__L_upper_2-torsion_function_for_chain_complexes}

 There is the notion of \emph{$L^2$-torsion}, which generalizes the classical notion of
 Reidemeister torsion to the $L^2$-setting. It was introduced topologically by
 L\"uck--Rothenberg~\cite{Lueck-Rothenberg(1991)} and analytically by
 Lott~\cite{Lott(1992a)} and Mathai~\cite{Mathai(1992)}. The identification of the
 analytical and topological $L^2$-torsion was carried out by
 Burghelea--Friedlander--Kappeler--McDonald~\cite{Burghelea-Friedlander-Kappeler-McDonald(1996a)}.
 For basic information about $L^2$-Betti numbers and $L^2$-torsion we refer for instance
 to~\cite{Kammeyer(2019), Lueck(2002)}.

 Versions, where a twisting with a (not necessarily surjective) homomorphism
 $\phi \colon G \to \IZ$ or more generally with a finite-dimensional representation is
 built in, were introduced and investigated by Dubois, Friedl, Liu, and L\"uck
 in~\cite{Dubois-Friedl-Lueck(2016), Friedl-Lueck(2019Thurston), Liu(2017), Lueck(2018)}.

 We review here the notion of \emph{$\phi$-twisted $L^2$-torsion} in the special case,
 where $G$ satisfies Condition~\ref{basic_setup}~\ref{basic_setup:G} below. This is enough
 for our purposes and is easy to state and handle. Namely, we do not have to deal with the
 difference between simple homotopy equivalences and homotopy equivalences and the
 question whether the chain complexes or spaces under consideration are $\phi$-twisted
 $\det$-$L^2$-acyclic.


 \subsection{Basic setup for \texorpdfstring{$\phi$}{phi}-twisted \texorpdfstring{$L^2$}{L2}-torsion function}%
 \label{subsec:Basic_setup_for_varphi-twisted_L_upper_2-torsion_function}

\begin{definition}[Farrell--Jones group]\label{def:Farrell-Jones_group}
  We call a group a \emph{Farrell--Jones group} if it satisfies the Full Farrell--Jones
  Conjecture.
\end{definition}

The statement of the Full Farrell--Jones Conjecture is quite complicated and involves
$L$-theory as well as $K$-theory, see~\cite[Conjecture~13.30 on page~387]{Lueck(2025book)}.
Its precise formulation is not relevant for this paper. For a status report and a list of
applications we refer to~\cite[Theorem~16.1 on page~481 and Theorem~13.65 on
page~405]{Lueck(2025book)}. Note that a hyperbolic group, a finite-dimensional
CAT(0)-group, a lattice in a path connected locally compact second countable Hausdorff
group, and the fundamental group of a (not necessarily compact or orientable) $3$-manifold
(with or without boundary) satisfy the Full Farrell--Jones Conjecture and that the Full
Farrell--Jones Conjecture has very good inheritance properties, e.g., it passes to
subgroups, overgroups of finite index, direct and free products, and directed colimits
(with arbitrary structure maps).

We will consider the following situation:

\begin{basicsetup}\label{basic_setup}\

  \begin{enumerate}

  \item\label{basic_setup:G} A finitely generated residually finite group $G$ which is a
    Farrell--Jones group.

  \item\label{basic_setup:X} A finite free $\IZ[G]$-chain complex $C_*$ which is
    \emph{$L^2$-acyclic}, i.e., the $L^2$-Betti numbers
    $b_n^{(2)}(C_*;\caln(G)) := b_n^{(2)}(L^2(G) \otimes_{\IZ[G]} C_*;\caln(G))$ vanish
    for $n \in \IZ_{\ge 0}$;

  \item\label{basic_setup:varphi} A (not necessarily surjective) group homomorphism
    $\phi \colon G \to \IZ$.

  \end{enumerate}
\end{basicsetup}

Note that Condition~\ref{basic_setup}~\ref{basic_setup:G} passes to
finitely generated subgroups, overgroups of finite index, and direct and free products.

We call two (not necessarily continuous) functions $\rho_0, \rho_1 \colon \IR_{>0} \to \IR$ \emph{equivalent},
in signs $\rho_0 \sim \rho_1$, if there exists $m \in \IZ$ such that
$\rho_0(t) - \rho_1(t) = m \cdot \ln(t)$ holds for $t \in \IR_{> 0}$. Let
\begin{equation}
  \cale := \{\overline{\rho} \mid \rho \colon \IR_{> 0} \to \IR\}
  \label{cale}
\end{equation}
be the set of equivalence classes $\overline{\rho}$ of  functions
$\rho \colon \IR_{>0} \to \IR$. Note that this becomes an abelian group by
\begin{eqnarray*}
  \overline{\rho_0} +\overline{\rho_1}
  &:=&
       \overline{\rho_0 + \rho_1};
  \\
  - \overline{\rho}
  & := &
         \overline{-\rho};
  \\
  0
  & = &
        \overline{0},
\end{eqnarray*}
where $\rho_0 + \rho_1$ sends $t$ to $\rho_0(t) + \rho_1(t)$ and $-\rho$ sends $t$ to
$-\rho(t)$.  We will call $[\rho]$ \emph{continuous} if one (and hence every)
representative $\rho$ is continuous.

Given the data above, one can define an element called \emph{reduced $\phi$-twisted
  $L^2$-torsion function}
\begin{equation}
  \overline{\rho}^{(2)}(C_*;\phi) \in \cale.
  \label{reduced_twisted_L_upper_torsion_function_chain}
\end{equation}

The construction of $\overline{\rho}^{(2)}(C_*;\phi)$ and the proof that it is
well-defined can be found in~\cite[Theorem~7.7~(3) and Remark~8.7]{Lueck(2018)}.  The
precise definition is not needed to understand this paper and can be taken as a black box.
We include a brief account after the next definition for the interested reader.

\begin{definition}[Degree of an equivalence class of an element in $\cale$]%
  \label{def:Degree_of_an_equivalence_class_of_an_element_in_cale}
  Let $\overline{\rho}$ be an equivalence class of functions $\IR_{> 0} \to \IR$. Let
  $\rho$ be a representative.  Assume that
  $\liminf_{t \to 0+} \frac{\rho(t)}{\ln(t)} \in \IR$ and
  $\limsup_{t \to \infty} \frac{\rho(t)}{\ln(t)} \in \IR$.

  Then define the \emph{degree at zero} and the \emph{degree at infinity} of $\rho$ to be
  the  real numbers
  \begin{eqnarray*}
    \deg_0(\rho) 
    & := & 
           \liminf_{t \to 0+} \frac{\rho(t)}{\ln(t)};
    \\
    \deg_{\infty}(\rho) 
    & := & 
           \limsup_{t \to \infty} \frac{\rho(t)}{\ln(t)}.
  \end{eqnarray*}
  Define the \emph{degree} of $\overline{\rho}$ to be the real number, which is
  independent of the choice of representative,
  \begin{eqnarray*}
    \deg(\overline{\rho}) 
    & := & 
           \deg_{\infty}(\rho)  -  \deg_0(\rho).
  \end{eqnarray*}
\end{definition}


The remainder of the section is not necessary to understand the rest
of the article.  However, the reader made find it helpful to see the
construction of the $\phi$-twisted $L^2$-torsion function and so we
include it below.

The following is taken from\cite[\S1]{Lueck(2018)} An $R$-module $M$
is \emph{based} if it is free and we have chosen an equivalence
classes of $R$-bases $[B_M]$ for $M$.  Here $B$ and $B'$ are
\emph{equivalent} if there exists a bijection $\sigma\colon B\to B'$
and $\sigma(b)=\pm b'$.  For $V$ a finite dimensional left
$G$-representation over $\IC$ we say that $V$ is a \emph{based finite
  dimensional $G$-representation}.  Define the following categories:
\begin{enumerate}
	\item $R$-$\mathbf{FBMod}$, based finite generated free $R$-modules with $R$-linear maps;
	\item $R$-$\mathbf{FBCC}$, chain complexes over $R$-$\mathbf{FBMod}$;
	\item $\caln(G)$-$\mathbf{FGHil}$, the additive $\IC$-category of finitely generated Hilbert $\caln(G)$-modules;
	\item $\caln(G)$-$\mathbf{FGHCC}$, chain complexes over $\caln(G)$-$\mathbf{FGHil}$.
\end{enumerate}
There are functors 
\begin{align} \label{functor:FBMod_to_FBMod}
	\eta^G_{V,[B_V]}\colon  \IC[G]\text{-}\mathbf{FBMod} &\to  \IC[G]\text{-}\mathbf{FBMod} \\ 
	(M,[B_M]) &\mapsto (M\otimes_\IC V, \{b\otimes v\mid b\in B_M, v\in V\}) \nonumber\\
	f &\mapsto f\otimes_\IC \id_V \nonumber
\end{align}
which extend to 
\begin{equation}
	\eta^G_{V,[B_V]}\colon  \IC[G]\text{-}\mathbf{FBCC} \to  \IC[G]\text{-}\mathbf{FBCC}
\end{equation}
in the obvious way.  There is also a functor
\begin{align}
 	\Lambda^G\colon \IC[G]\text{-}\mathbf{FBMod} &\to \caln(G)\text{-}\mathbf{FGHil} \\
 	(M,[B_M]) &\mapsto L^2(G)\otimes_{\IC[G]}M \nonumber\\
 	f &\mapsto \id_{L^2(G)}\otimes_\IC f \nonumber 	
\end{align}
where we choose a $B_M$ such that the map $\bigoplus_{b\in B_M} L^2(G) \to L^2(G)\otimes_{\IC[G]} M$ given by $(x_b)\mapsto \sum x_b \otimes b$ is an isometry.  This extends to a functor
\begin{equation}
 	\Lambda^G\colon  \IC[G]\text{-}\mathbf{FBCC} \to  \caln(G)\text{-}\mathbf{FBCC}
\end{equation}
in the obvious way.

\begin{remark}
  In each of the above functors we will often suppress the notation when it is clear from
  context, e.g. $\eta^G_{V,[B_V]}$ may be denoted by $\eta$.
\end{remark}

\begin{example}
  Let $G$ be a group and $\phi\colon G\to \IR$ a homomorphism.  Let $\IC_t$ be the based
  one-dimensional $\IR$-representation where $r\in\IR$ acts by $t^r$ and $\IC_t$ has basis
  $\{t\}$.  Denote by $\phi^\ast\IC_t$ the based one-dimensional $G$-representation
  obtained from $\IC_t$.  The functor~\eqref{functor:FBMod_to_FBMod} behaves as follows:
  Let $A\in M_{r,s}(\IZ[G])$, we have $r_A\colon\IZ[G]^r\to \IZ[G]^s$ (by right
  multiplication) is sent under $\eta_{\phi^\ast\IC_t}$ to
  $r_{\hat\phi_\ast(A)}\colon \IC[G]^r\to \IC[G]^s$, where the $(r,s)$-matrix
  $\hat\phi_\ast(A)\in M_{r,s}(\IC[G])$ is obtained elementwise by
  \begin{align*}
    \hat\phi\colon \IC[G] &\to\IC[G] \\
    \sum_{g\in G}\lambda_g g &\mapsto \sum_{g\in G}t^{\phi(g)}\lambda_g g.
  \end{align*}
\end{example}

Let $G$ be a group and let $\phi\colon G\to \IR$ be a group homomorphism.  Let
$(C_\ast,[B_{C_\ast}])$ be a based $L^2$-acyclic chain complex of free $\IZ[G]$-modules.
We say $C_\ast$ is \emph{$\phi$-twisted determinant class} or \emph{$\phi$-twisted
  $L^2$-acyclic} if the finite $\caln(G)$-chain complex
$\Lambda\circ\eta_{\phi^\ast\IC_t}$ is of determinant class or $L^2$-acyclic.

Note that if $G$ is residually finite, then $C_\ast$ being $L^2$-acyclic automatically
implies that $C_\ast$ is both $\phi$-twisted $L^2$-acyclic and $\phi$-twisted determinant
class.

Suppose that $C_\ast$ is $\phi$-twisted determinant class.  We define the $\phi$-twisted
$L^2$-torsion function by
\begin{align*}
\rho^{(2)}(C_\ast,[B_{C_\ast}];\phi)\colon \IR_{>0} &\to \IR  \\
	t &\mapsto \rho^{(2)}(\Lambda\circ\eta_{\phi^\ast\IC_t}(C_\ast,[B_{C_\ast}])),
\end{align*}
that is we send $t\in\IR_{>0}$ to the $L^2$-torsion of $C_\ast$ twisted with $\phi^\ast\IC_t$.

The reduced $\phi$-twisted $L^2$-torsion does not depend on the basis refinement
$[B_{C_\ast}]$ so we will usually suppress this from the notation.


\subsection{Basic properties of \texorpdfstring{$\phi$}{phi}-twisted \texorpdfstring{$L^2$}{L2}-torsion function}%
\label{subsec:Basic_properties_of_varphi-twisted_L_upper_2-torsion_function}

\begin{theorem}\label{the:basic_properties_of_overline(rho)(C_ast;varphi)}
  Let $G$ be a group satisfying Condition~\ref{basic_setup}~\ref{basic_setup:G}, i.e., $G$
  is finitely generated and residually finite and satisfies the Full Farrell--Jones
  Conjecture. Let $\phi \colon G \to \IZ$ be a group homomorphism.
  
  \begin{enumerate}
  \item\label{the:basic_properties_of_overline(rho)(C_ast;varphi)_homotopy_invariance} Let
    $C_*$ and $D_*$ be finite free $\IZ[G]$-chain complexes which are $\IZ[G]$-chain
    homotopy equivalent.  Suppose that $C_*$ or $D_*$ is $L^2$-acyclic.

    Then both $C_*$ and $D_*$ are $L^2$-acyclic and we get
    \[
      \overline{\rho}^{(2)}(C_*;\phi) = \overline{\rho}^{(2)}(D_*;\phi);
    \]
  \item\label{the:basic_properties_of_overline(rho)(C_ast;varphi):Additivity} Let
    $0 \to C_* \to D_* \to E_* \to 0$ be an exact sequence of finite free $\IZ[G]$-chain
    complexes.  Suppose that two of the $\IZ[G]$-chain complexes $C_*$, $D_*$, and $E_*$
    are $L^2$-acyclic.

    Then all three are $L^2$-acyclic and we get
    \[
      \overline{\rho}^{(2)}(D_*;\phi) = \overline{\rho}^{(2)}(C_*;\phi) +
      \overline{\rho}^{(2)}(E_*;\phi);
    \]

  \item\label{the:basic_properties_of_overline(rho)(C_ast;varphi):product_formula} Let
    $C_*$ be a finite free $\IZ[G]$-chain complex which is $L^2$-acyclic.  Let $H$ be a
    group satisfying Condition~\ref{basic_setup}~\ref{basic_setup:G} and let $D_*$ be a
    finite free $\IZ[H]$-chain complex. Let $\pr_G \colon G \times H \to G$ be the
    projection and $\chi_{\IZ[H]}(D_*) \in \IZ$ be the Euler characteristic of $D_*$.

    Then $G \times H$ satisfies Condition~\ref{basic_setup}~\ref{basic_setup:G},
    $C_* \otimes_{\IZ} D_*$ is a finite free $\IZ[G \times H]$-chain complex which is
    $L^2$-acyclic, and we get
    \[
      \overline{\rho}^{(2)}(C_* \otimes_{\IZ} D_*;\phi \circ \pr_G) = \chi_{\IZ[H]}(D_*)
      \cdot \overline{\rho}^{(2)}(C_*;\phi);
    \]

  \item\label{the:basic_properties_of_overline(rho)(C_ast;varphi):induction} Let
    $i \colon H \to G$ be the inclusion of the finitely generated subgroup $H$ of $G$.
    Let $C_*$ be a finite free $\IZ[H]$-chain complex which is $L^2$-acyclic.

    Then $H$ satisfies Condition~\ref{basic_setup}~\ref{basic_setup:G},
    $i_*C_* = \IZ[G] \otimes_{\IZ[H]} C_*$ is a finite free $\IZ[G]$-chain complex which
    is $L^2$-acyclic, and we get
    \[
      \overline{\rho}^{(2)}(C_*;\phi \circ i) = \overline{\rho}^{(2)}(i_*C_*;\phi);
    \]
    
  \item\label{the:basic_properties_of_overline(rho)(C_ast;varphi);restriction} Let
    $i \colon H \to G$ be the inclusion of the subgroup $H$ of $G$ with finite index
    $[G:H]$.  Let $C_*$ be a finite free $\IZ[G]$-chain complex which is $L^2$-acyclic.

    Then $H$ satisfies Condition~\ref{basic_setup}~\ref{basic_setup:G}, the $\IZ[H]$-chain
    complex $i^*C_*$ obtained from $C_*$ by restriction with $i$ is a finite free
    $\IZ[H]$-chain complex which is $L^2$-acyclic, and we get
    \[
      \overline{\rho}^{(2)}(i^*C_*;\phi \circ i) = [G:H] \cdot
      \overline{\rho}^{(2)}(C_*;\phi);
    \]

  \item\label{the:basic_properties_of_overline(rho)(C_ast;varphi):scaling} Let $C_*$ be a
    finite free $\IZ[G]$-chain complex which is $L^2$-acyclic. Consider
    $\phi \colon G \to \IZ$.  For $m \in \IZ$ let $\phi_m \colon G \to \IZ$ be the
    composite of $\phi$ with $m \cdot \id_{\IZ}$.  Let $\rho$ be a representative of
    $\overline{\rho}^{(2)}(C_*;\phi)$. Define $\rho_m(t) := \rho(t^m)$.

    Then $\rho_m(t)$ is a representative of $\overline{\rho}^{(2)}(C_*;\phi_m)$;

  \item\label{the:basic_properties_of_overline(rho)(C_ast;varphi):value_at_one} If $\rho$
    is a representative of $\overline{\rho}^{(2)}(C_*;\phi)$, then $\rho(1)$ agrees with
    the $L^2$-torsion
    $\rho^{(2)}(C_*;\caln(G)) := \rho^{(2)}(L^2(G) \otimes_{\IZ[G]} C_*;\caln(G))$;

  \item\label{the:basic_properties_of_overline(rho)(C_ast;varphi):trivial_twisting} If
    $\phi \colon G \to \IZ$ is trivial, then the constant function with value the
    $L^2$-torsion $\rho^{(2)}(C_*;\caln(G))$ is a representative of
    $\overline{\rho}^{(2)}(C_*;\phi)$;

  \item\label{the:basic_properties_of_overline(rho)(C_ast;varphi):continuous}
    $\overline{\rho}^{(2)}(C_*;\phi)$ is continuous;

  \item\label{the:basic_properties_of_overline(rho)(C_ast;varphi):degree} The degree of
    $\overline{\rho}^{(2)}(C_*;\phi)$ is defined.

  \end{enumerate}
\end{theorem}
\begin{proof}
  Theorem~\ref{the:basic_properties_of_overline(rho)(C_ast;varphi)} is just the chain
  complex version of~\cite[Theorem~8.6 and Theorem~10.2]{Lueck(2018)}
  taking~\cite[Theorem~8.3~(1) and (2)]{Lueck(2018)} into account.
\end{proof}


\subsection{{Review on \texorpdfstring{$\phi$}{phi}-twisted \texorpdfstring{$L^2$}{L2}-torsion function for \texorpdfstring{$CW$}{CW}-complexes}}%
\label{sec:Review_on_varphi-twisted__L_upper_2-torsion_function_for_CW-complexes}

Let $G$ be a group satisfying Condition~\ref{basic_setup}~\ref{basic_setup:G}, i.e., $G$
is finitely generated and residually finite and satisfies the Full Farrell--Jones
Conjecture.  Let $\phi \colon G \to \IZ$ be a (not necessarily surjective) homomorphism.
Let $Y$ be a finite free $G$-$CW$-complex which is \emph{$L^2$-acyclic} in the sense that
all its $L^2$-Betti numbers $b_n^{(2)}(Y;\caln(G)) = b_n^{(2)}(C_*(Y);\caln(G))$ vanish.
Then we can define its \emph{reduced $\phi$-twisted $L^2$-torsion function}
\begin{equation}
  \overline{\rho}^{(2)}(Y;\phi) \in \cale
  \label{reduced_twisted_L_upper_torsion_function_for_CW}
\end{equation}
to be  $\overline{\rho}^{(2)}(Y;\phi) := \overline{\rho}^{(2)}(C_*(Y);\phi)$ for the reduced
$\phi$-twisted $L^2$-torsion function $\overline{\rho}^{(2)}(C_*(Y);\phi)$ of its cellular chain
$\IZ[G]$-complex $C_*(Y)$, see~\eqref{reduced_twisted_L_upper_torsion_function_chain}.

If $X$ is a connected finite $CW$-complex with fundamental group $G = \pi_1(X)$, we can
take $Y$ to be the universal covering $\widetilde{X}$ and we write
\begin{equation}
  \overline{\rho}^{(2)}(\widetilde{X};\phi) \in \cale.
  \label{reduced_twisted_L_upper_torsion_function_for_universal_covering}
\end{equation}

Consider a (not necessarily connected) finite $CW$-complex $X$ and an element
$\phi \in H^1(X;\IZ)$. We will abuse notation and will say that $\pi_1(X)$ satisfies
Condition~\ref{basic_setup}~\ref{basic_setup:G} if $\pi_1(C)$ satisfies
Condition~\ref{basic_setup}~\ref{basic_setup:G} for every path component $C \in \pi_0(X)$,
and that $\widetilde{X}$ is \emph{$L^2$-acyclic} if $\widetilde{C}$ is $L^2$-acyclic for
every $C \in \pi_0(X)$. If $\pi_1(X)$ satisfies
Condition~\ref{basic_setup}~\ref{basic_setup:G} and $\widetilde{X}$ is $L^2$-acyclic, then
we can define

\begin{equation}
  \overline{\rho}^{(2)}(\widetilde{X};\phi)
  = \sum_{C \in \pi_0(X)} \overline{\rho}^{(2)}(\widetilde{C};\phi|_C) \in \cale,
  \label{reduced_twisted_L_upper_torsion_function_for_universal_covering_not_connected}
\end{equation}
where $\phi|_C \in H^1(C;\IZ)$ is obtained from $\phi$ by restriction to the component $C$
and we think of $\phi|_C$ as a group homomorphism $\phi|_C \colon \pi_1(C) \to \IZ$.  Note
that~\eqref{reduced_twisted_L_upper_torsion_function_for_universal_covering_not_connected}
reduces to~\eqref{reduced_twisted_L_upper_torsion_function_for_universal_covering} if $X$
is connected.


\subsection{Basic properties of the reduced \texorpdfstring{$L^2$}{L2}-torsion function for universal
  coverings}%
\label{subsec:Basic_properties_of_the_reduced_L2-torsion_function_for_universal_coverings}

\begin{theorem}[Basic properties of the reduced $L^2$-torsion function for universal
  coverings]%
  \label{the:Basic_properties_of_the_reduced_L2-torsion_function_for_universal_coverings}
  \

  \begin{enumerate}

  \item\label{the:Basic_properties_of_the_reduced_L2-torsion_function_for_universal_coverings:homotopy_invariance}
    \emph{Homotopy invariance}\\
    Let $f \colon X \to Y$ be a homotopy equivalence of finite $CW$-complexes.  Consider
    $\phi_X \in H^1(X;\IZ)$ and $\phi_Y \in H^1(Y;\IZ)$.  Suppose that $\pi_1(X)$
    satisfies Condition~\ref{basic_setup}~\ref{basic_setup:G}, that $\widetilde{X}$ is
    $L^2$-acyclic, and that $\phi_X = H^1(f;\IZ)(\phi_Y)$ holds.

    Then $\pi_1(Y)$ satisfies Condition~\ref{basic_setup}~\ref{basic_setup:G} and
    $\widetilde{Y}$ is $L^2$-acyclic, and we get in $\cale$
    \[
      \overline{\rho}^{(2)}(\widetilde{X};\phi_X) =
      \overline{\rho}^{(2)}(\widetilde{Y};\phi_Y);
    \]

  \item\label{the:Basic_properties_of_the_reduced_L2-torsion_function_for_universal_coverings:sum_formula}
    \emph{Sum formula}\\
    Consider a cellular pushout of finite $CW$-complexes
    \[
      \xymatrix{ X_0 \ar[r]^{i_1} \ar[rd]^{j_0} \ar[d]_{i_2} & X_1 \ar[d]^{j_1}
        \\
        X_2 \ar[r]_{j_2} & X }
    \]
    where $i_1$ is cellular, $i_2$ is an inclusion of $CW$-complexes, and $X$ has the
    obvious $CW$-structure coming from the ones on $X_0$, $X_1$, and $X_2$. Suppose that
    for $i =0,1,2$ the map $j_i$ is $\pi_1$-injective, i.e., for any choice of base point
    $x_i \in X_i$ the induced map
    $\pi_1(j_i,x_i) \colon \pi_1(X_i,x_i) \to \pi_1(X,j_i(x_i))$ is injective.  Suppose we
    are given elements $\phi_i \in H^1(X_i;\IZ)$ and $\phi \in H^1(X;\IZ)$ such that
    $j_i^*(\phi) = \phi_i$ holds for $i = 0,1,2$.  Assume that $\widetilde{X_i}$ is
    $L^2$-acyclic for $i = 0,1,2$ and that $\pi_1(X)$ satisfies
    Condition~\ref{basic_setup}~\ref{basic_setup:G}.

    Then $\pi_1(X_i)$ satisfies Condition~\ref{basic_setup}~\ref{basic_setup:G} for
    $i = 0,1,2$ and $\widetilde{X}$ is $L^2$-acyclic and we get
    \[
      \overline{\rho}^{(2)}(\widetilde{X};\phi) =
      \overline{\rho}^{(2)}(\widetilde{X_1};\phi_1) +
      \overline{\rho}^{(2)}(\widetilde{X_2};\phi_2) -
      \overline{\rho}^{(2)}(\widetilde{X_0};\phi_0);
    \]

  \item\label{the:Basic_properties_of_the_reduced_L2-torsion_function_for_universal_coverings:product_formula}
    \emph{Product formula}\\
    Let $X$ and $Y$ be finite $CW$-complexes. Suppose that $Y$ is connected and that
    $\pi_1(X)$ and $\pi_1(Y)$ satisfy Condition~\ref{basic_setup}~\ref{basic_setup:G}.
    Consider an element $\phi \in H^1(X \times Y;\IZ)$. Define $\phi_X$ to be the image of
    $\phi$ under the map $H^1(X \times Y;\IZ) \to H^1(X;\IZ)$ induced by the inclusion
    $X \to X \times Y, \; x \mapsto (x,y)$ for any choice of base point $y \in Y$.
    Suppose that $\widetilde{X}$ is $L^2$-acyclic.

    Then $\widetilde{X \times Y}$ is $L^2$-acyclic, $\pi_1(X \times Y)$ satisfies
    Condition~\ref{basic_setup}~\ref{basic_setup:G}, and we get
    \[
      \overline{\rho}^{(2)}(\widetilde{X \times Y}; \phi) = \chi(Y) \cdot
      \overline{\rho}^{(2)}(\widetilde{X};\phi_X);
    \]

  \item\label{the:Basic_properties_of_the_reduced_L2-torsion_function_for_universal_coverings:Poincare_duality}
    \emph{Poincar\'e duality}\\
    Let $X$ be a finite $n$-dimensional Poincar\'e complex, e.g., a closed manifold of
    dimension $n$. Consider $\phi \in H^1(X;\IZ)$.  Suppose that $X$ is $L^2$-acyclic and
    $\pi_1(X)$ satisfies Condition~\ref{basic_setup}~\ref{basic_setup:G}.
  
    Then
    \[
      \overline{\rho}^{(2)}(\widetilde{X};\phi)(t) = (-1)^{n+1} \cdot
      \overline{\rho}^{(2)}(\widetilde{X};\phi)(t^{-1});
    \]

  \item\label{the:Basic_properties_of_the_reduced_L2-torsion_function_for_universal_coverings:(Multiplicativity)}
    \emph{Finite coverings}\\
    Let $p \colon X \to Y$ be a $d$-sheeted covering of finite $CW$-complexes for some
    natural number $d$.  Let $\phi_Y \in H^1(Y;\IZ)$ and $\phi_X \in H^1(X;\IZ)$ be
    elements satisfying $H^1(p;\IZ)(\phi_Y) = \phi_X$.

    Then $\widetilde{X}$ is $L^2$-acyclic and $\pi_1(X)$ satisfies
    Condition~\ref{basic_setup}~\ref{basic_setup:G}, if and only if $\widetilde{Y}$ is
    $L^2$-acyclic and $\pi_1(Y)$ satisfies
    Condition~\ref{basic_setup}~\ref{basic_setup:G}, and in this case we get
    \[
      \overline{\rho}^{(2)}(\widetilde{X};\phi_X) = d \cdot
      \overline{\rho}^{(2)}(\widetilde{Y};\phi_Y);
    \]

  \item\label{the:Basic_properties_of_the_reduced_L2-torsion_function_for_universal_coverings:scaling}
    \emph{Scaling $\phi$}
    \\
    Let $m \in \IZ$,  and let $X$ be a finite $CW$-complex $X$ such that $\pi_1(X)$ satisfies
    Condition~\ref{basic_setup}~\ref{basic_setup:G} and $\widetilde{X}$ is $L^2$-acyclic,
    and $\phi \in H^1(X;\IZ)$.  Let $\rho$ be a representative of
    $\overline{\rho}^{(2)}(\widetilde{X};\phi)$. Put $\rho_m(t) = \rho(t^m)$.

    Then $\rho_m$ is a representative of
    $\overline{\rho}^{(2)}(\widetilde{X}; m\cdot \phi)$;

  \item\label{the:Basic_properties_of_the_reduced_L2-torsion_function_for_universal_coveringsvalue_at_one}
    Let $X$ be a finite $CW$-complex $X$ such that $\pi_1(X)$ satisfies
    Condition~\ref{basic_setup}~\ref{basic_setup:G} and $\widetilde{X}$ is $L^2$-acyclic,
    and $\phi \in H^1(X;\IZ)$.  Let $\rho$ be  a representative of
    $\overline{\rho}^{(2)}(\widetilde{X};\phi)$.

    Then $\rho(1)$ agrees with the $L^2$-torsion
    $\rho^{(2)}(\widetilde{X};\caln(\pi_1(X)))$;

  \item\label{tthe:Basic_properties_of_the_reduced_L2-torsion_function_for_universal_coverings:trivial_twisting}
    Let $X$ be a finite $CW$-complex $X$ such that $\pi_1(X)$ satisfies
    Condition~\ref{basic_setup}~\ref{basic_setup:G} and $\widetilde{X}$ is $L^2$-acyclic.
    Suppose that $\phi \colon \pi_1(X) \to \IZ$ is trivial.

    Then the constant function with value $\rho^{(2)}(\widetilde{X};\caln(\pi_1(X)))$ is a
    representative of $\overline{\rho}^{(2)}(\widetilde{X};\phi)$;

  \item\label{the:Basic_properties_of_the_reduced_L2-torsion_function_for_universal_coverings:continuous}
    Let $X$ be a finite $CW$-complex $X$ such that $\pi_1(X)$ satisfies
    Condition~\ref{basic_setup}~\ref{basic_setup:G} and $\widetilde{X}$ is $L^2$-acyclic,
    and $\phi \in H^1(X;\IZ)$.

    Then $\overline{\rho}^{(2)}(\widetilde{X};\phi)$ is continuous;

  \item\label{the:Basic_properties_of_the_reduced_L2-torsion_function_for_universal_coverings:degree}
    Let $X$ be a finite $CW$-complex $X$ such that $\pi_1(X)$ satisfies
    Condition~\ref{basic_setup}~\ref{basic_setup:G} and $\widetilde{X}$ is $L^2$-acyclic,
    and $\phi \in H^1(X;\IZ)$.

    Then the degree of $\overline{\rho}^{(2)}(\widetilde{X};\phi)$ is defined.

  \end{enumerate}
\end{theorem}
\begin{proof}
  Theorem~\ref{the:Basic_properties_of_the_reduced_L2-torsion_function_for_universal_coverings}
  follows from~\cite[Theorem~8.15 and Theorem~10.2]{Lueck(2018)}
  taking~\cite[Theorem~8.3~(1) and (2)]{Lueck(2018)} into account.
\end{proof}


\subsection{Fibration formula}\label{subsec:fibration_formula}

We also have the following results generalizing the product formula appearing in
Theorem~\ref{the:Basic_properties_of_the_reduced_L2-torsion_function_for_universal_coverings}~%
\ref{the:Basic_properties_of_the_reduced_L2-torsion_function_for_universal_coverings:product_formula}.

\begin{theorem}[Fibrations]\label{the:fibrations_rho}
  Let $F \xrightarrow{i} E \xrightarrow{p} B$ be a fibration such that $F$ and $B$ are
  connected finite $CW$-complexes and $\pi_1(i) \colon \pi_1(F) \to \pi_1(E)$ is
  injective. Consider $\phi \in H^1(E;\IZ)$.  Let $i^*\phi \in H^1(F;\IZ)$ be its pullback
  with $i$.  Suppose that $\widetilde{F}$ is $L^2$-acyclic and $\pi_1(E)$ satisfies
  Condition~\ref{basic_setup}~\ref{basic_setup:G}.

  Then $\pi_1(F)$ satisfies Condition~\ref{basic_setup}~\ref{basic_setup:G}, $E$ is up to
  homotopy a finite $CW$-complex, $\widetilde{E}$ is $L^2$-acyclic, and we get
  \[
    \overline{\rho}^{(2)}(\widetilde{E};\phi) = \chi(B) \cdot
    \overline{\rho}^{(2)}(\widetilde{F};i^*\phi).
  \]
\end{theorem}
\begin{proof}
  See~\cite[Theorem~8.16]{Lueck(2018)}.
\end{proof}

The next result follows from
Theorem~\ref{the:Basic_properties_of_the_reduced_L2-torsion_function_for_universal_coverings}~%
\ref{the:Basic_properties_of_the_reduced_L2-torsion_function_for_universal_coverings:scaling},
Theorem~\ref{the:fibrations_rho} and~\cite[(3.23) on page~136]{Lueck(2002)}.

\begin{corollary}\label{cor:S_uuper_1-fibrations}
  Let $S^1 \xrightarrow{i} E \xrightarrow{p} B$ be a fibration of connected finite
  $CW$-complexes such that $\pi_1(i) \colon \pi_1(S^1) \to \pi_1(E)$ is injective and
  $\pi_1(E)$ satisfies Condition~\ref{basic_setup}~\ref{basic_setup:G}.  Consider
  $\phi \in H^1(E;\IZ)$.  Let $i^*\phi \in H^1(S^1;\IZ)$ be its pullback with $i$.  Let
  $m \in \IZ$ be the integer for which $i^*\phi = m \cdot [S^1]$ holds for the standard
  generator $[S^1]$ of the infinite cyclic group $H^1(S^1;\IZ)$.

  Then $E$ is up to homotopy a finite $CW$-complex, $\widetilde{E}$ is $L^2$-acyclic, and
  the function
  \[
    \rho(t) = \begin{cases} 0 & t \le 1;
      \\
      \chi(B) \cdot |m| \cdot \ln(t) & t \ge 1,
    \end{cases}
  \]
  is a representative for $\overline{\rho}^{(2)}(\widetilde{E};\phi)$.
\end{corollary}


\typeout{--------------------- Section 3: Spectral radii and continuity -------------------}

\section{Spectral radii and continuity}%
\label{sec:Twisted_L2-torsion_function_and_mapping_tori_for_chain_complexes}

In this section we investigate the $\phi$-twisted $L^2$-torsion function of mapping tori
of chain maps.


\subsection{\texorpdfstring{$L^1$}{L1}-norm and spectral
  radius}\label{subsec:L_upper_1-norm_and_spectra_radius}

For an element $x = \sum_{g \in g} r_g \cdot g$ in the complex group ring $\IC G$ of the
group $G$ define its \emph{$L^1$-norm} to be
\begin{equation}
  ||x||_1 := \sum_{g \in G} |r_g| \in \IR_{\ge 0}.
  \label{L_1_norm_for_elements_in_CG}
\end{equation}

Let $A \in \Mpar{m}{n}{\IC[G]}$ be an $(m,n)$ matrix over $\IC[G]$. So $A$ looks like
\[ A =
  \begin{pmatrix}
    a_{1,1} & a_{1,2} & a_{1,3} & \cdots & a_{1 ,n-2} & a_{1,n-1} & a_{1,n}
    \\
    a_{2,1} & a_{2,2} & a_{2,3} & \cdots & a_{2,n-2} & a_{2,n-1} & a_{2,n}
    \\
    a_{3,1} & a_{3,2} & a_{3,3} & \cdots & a_{3,n-2} & a_{3,n-1} & a_{3,n}
    \\
    \vdots & \vdots & \vdots & \ddots &\vdots & \vdots & \vdots
    \\
    a_{m-1,1} & a_{m-1,2} & a_{m-1,3} & \cdots & a_{m-1,n-2} & a_{m-1,n-1} & a_{m-1,n}
    \\
    a_{m,1} & a_{m,2} & a_{m,3} & \cdots & a_{m,n-2} & a_{m,n-1} & a_{m,n}
  \end{pmatrix}.
\]

We define for a matrix $A \in \Mpar{m}{n}{\IC[G]}$ its \emph{$L^1$-norm}
\begin{equation}
  ||A||_1 := \max\left\{\left.\sum_{j = 1}^n||a_{i,j}||_1 \;\right|\; i \in \{1,2, \ldots, m\}\right\}
  \label{||A||_1_for__(m,n)_matrices}
\end{equation}
and its \emph{operator norm}
\begin{equation}
  ||A||_{\infty} := ||r_A^{(2)}||_{\infty}
  \label{||A||_infty_for_(m,n)_matrices}
\end{equation}
where $||r_A^{(2)}||_{\infty}$ is the operator norm of the bounded operator
$r_A^{(2)} \colon L^2(G)^m \to L^2(G)^n$ given by right multiplication with $A$.

The elementary proof of the next lemma is left to the reader.

  \begin{lemma}\label{lem:properties_of_the_L_uper_1-norm_for_matrices}
    Consider $A,B \in \Mpar{m}{n}{\IC[G]}$ and $\lambda \in \IC$.
    \begin{enumerate}
    \item\label{lem:properties_of_the_L_uper_1-norm_for_matrices:||A||_is_0} We have
      $||A||_1 = 0 \Longleftrightarrow ||A||_{\infty} = 0 \Longleftrightarrow A = 0$;
    \item\label{lem:properties_of_the_L_uper_1-norm_for_matrices:||I_n||_is_n} We have
      $||I_n||_1 = ||I_n||_{\infty} = 1$ for the unit $(n,n)$-matrix $I_n$;
    \item\label{lem:properties_of_the_L_uper_1-norm_for_matrices:scalar} We have
      $||\lambda \cdot A||_1 = |\lambda| \cdot ||A||_1$ and
      $||\lambda \cdot A||_{\infty} = |\lambda| \cdot ||A||_{\infty}$;
    \item\label{lem:properties_of_the_L_uper_1-norm_for_matrices:triangle} We have
      $||A + B||_1 \le ||A||_1 + ||B||_1$ and
      $||A + B||_{\infty} \le ||A||_{\infty} + ||B||_{\infty}$;
    \item\label{lem:properties_of_the_L_uper_1-norm_for_matrices:Banach} We have
      $||A \cdot C||_1 \le ||A||_1 \cdot ||C||_1$ and
      $||A \cdot C||_{\infty} \le ||A||_{\infty} \cdot ||C||_{\infty}$ for
      $A\in \Mpar{m}{n}{\IC[G]}$ and $C \in \Mpar{n}{o}{\IC[G]}$;
    \item\label{lem:properties_of_the_L_uper_1-norm_for__matrices:inequality} We have
      $||A||_{\infty} \le ||A||_1$.
      
    \end{enumerate}
  \end{lemma}

  \begin{definition}[Spectral radius]\label{def:spectral-radius}
    We define the \emph{spectral radius with respect to the operator norm} of
    $A \in \Mpar{n}{n}{\IC[G]}$ to be
    \[
      h(A) := \lim_{k \to \infty} ||A^k||_{\infty}^{1/k} \in \IR_{\ge 0}.
    \]
  \end{definition}
  Note that the limit exists by Fekete's Subadditivity Lemma and
  Lemma~\ref{lem:properties_of_the_L_uper_1-norm_for_matrices}%
~\ref{lem:properties_of_the_L_uper_1-norm_for_matrices:Banach} which together imply
  \[
    \lim_{k \to \infty} \frac{\ln(||A^k||_{\infty})}{k} =
    \inf\left\{\frac{\ln(||A^k||_{\infty})}{k} \mid k \in \IZ_{\ge 1}\right\}.
  \]

  \begin{lemma}\label{lem:basic_properties_of_spectral_radius} \
    \begin{enumerate}
    \item\label{lem:basic_properties_of_spectral_radius:inequality} We have
      $h(A) \le ||A||_{\infty} \le ||A||_1$ for $A \in \Mpar{n}{n}{\IC[G]}$;
    \item\label{lem:basic_properties_of_spectral_radius:conjugation} Consider
      $A \in \Mpar{n}{n}{\IC[G]}$ and $U \in \GLpar{n}{\IC[G]}$. Then we get
      \[
        h(UAU^{-1}) = h(A);
      \]
    \item\label{lem:basic_properties_of_spectral_radius:von_Neumann_series} If $h(A) < 1$,
      then the operator $r_{I_n -A}^{(2)} \colon L^2(G)^n \to L^2(G)^n$ is invertible with
      inverse $\sum_{i = 0}^{\infty} r^{(2)}_{A^i}$, where the sum is to be understood
      with respect to the operator norm.
    \end{enumerate}
  \end{lemma}
  \begin{proof}~\ref{lem:basic_properties_of_spectral_radius:inequality} We get
    $||A^k||_{\infty} \le ||A||_{\infty}^k$ from
    Lemma~\ref{lem:properties_of_the_L_uper_1-norm_for_matrices}~%
\ref{lem:properties_of_the_L_uper_1-norm_for_matrices:Banach}.  This implies
    $||A^k||_{\infty}^{1/k} \le ||A||_{\infty}$ and hence $h(A) \le ||A||_{\infty}$.  We
    get $||A||_{\infty} \le ||A||_1$ from
    Lemma~\ref{lem:properties_of_the_L_uper_1-norm_for_matrices}~%
\ref{lem:properties_of_the_L_uper_1-norm_for__matrices:inequality}.
    \\[1mm]~\ref{lem:basic_properties_of_spectral_radius:conjugation} We get from
    Lemma~\ref{lem:properties_of_the_L_uper_1-norm_for_matrices}~%
\ref{lem:properties_of_the_L_uper_1-norm_for_matrices:Banach}
    \[
      ||(UAU^{-1})^k||_{\infty}^{1/k} = ||UA^kU^{-1}||_{\infty}^{1/k} \le (||U||_{\infty}
      \cdot ||U^{-1}||_{\infty})^{1/k} \cdot ||A^k||_{\infty}^{1/k}.
    \]
    Since $(||U||_{\infty} \cdot ||U^{-1}||_{\infty}) > 0$ holds by
    Lemma~\ref{lem:properties_of_the_L_uper_1-norm_for_matrices}~%
\ref{lem:properties_of_the_L_uper_1-norm_for_matrices:||A||_is_0}, we get
    \[
      \lim_{k \to \infty} (||U||_{\infty} \cdot ||U^{-1}||_{\infty})^{1/k} = 1.
    \]
    This implies $h(UAU^{-1}) \le h(A)$. Since $U^{-1}(UAU^{-1})(U^{-1})^{-1} = A$ holds,
    we also get $h(A) \le h(UAU^{-1})$. Hence we obtain $h(A) = h(UAU^{-1})$.
    \\[1mm]~\ref{lem:basic_properties_of_spectral_radius:von_Neumann_series} See for
    instance~\cite[Theorem~3.3.3 on page~204]{Kadison-Ringrose(1997I)}.
  \end{proof}

  Note that the definition of $h(f)$ makes also sense for an endomorphism
  $f \colon M \to M$ of a finitely generated free $\IZ[G]$-module because of
  Lemma~\ref{lem:basic_properties_of_spectral_radius}~\ref{lem:basic_properties_of_spectral_radius:conjugation}.
  Namely, choose a $\IZ[G]$-isomorphism $\mu \colon M \to \IZ[G]^a$ for some
  $a \in \IZ_{\ge 0}$ and a matrix $A \in \Mpar{a}{a}{\IZ[G]}$ such that
  $f = \mu^{-1} \circ r_A \circ \mu$ holds and define
  \begin{equation}
    h(f) = h(A).
    \label{h(f)}
  \end{equation}


  \subsection{Continuity of the Fuglede--Kadison determinant}%
  \label{subsec:Continuity_of_the_Fuglede--Kadison_determinant}

  For an introduction to the \emph{Fuglede--Kadison determinant}
  ${\det}^{(2)}_{\caln(G)}(T)$ of a bounded $G$-equivariant operator
  $T \colon L^2(G)^m \to L^2(G)^n$ we refer for instance
  to~\cite[Section~3.2]{Lueck(2002)}. It is possible for the reader to deal with this
  notion as a black box assigning to $T$ a non-negative real number without knowing the
  precise definition and properties.  If $G$ is trivial and $T$ is injective, it is the
  square root of the classical determinant of $T^*T$.

\begin{lemma}\label{lem:continuity_of_det_on_GL_n(caln(G))}
  The map
  \[
    \GLpar{n}{\caln(G)} \to \IR, \quad A \mapsto {\det}_{\caln(G)}\bigl(r_A \colon
    L^2(G)^n \to L^2(G)^n\bigr)
  \]
  is continuous on the open subset 
  \[\GLpar{n}{\caln(G)} \subseteq  \Mpar{n}{n}{\caln(G)} =\calb(L^2(G)^n,L^2(G)^n)^G,\]
  where we equip $\calb(L^2(G)^n,L^2(G)^n)^G$ with the  operator norm.
\end{lemma}
\begin{proof}
  See~\cite[Theorem~1.10~(d)]{Carey-Farber-Mathai(1997)} or~\cite[Theorem~1
  (3)]{Fuglede-Kadison(1952)}.
\end{proof}

We also mention the following result~\cite[Theorem~1.10~(e)]{Carey-Farber-Mathai(1997)}.

\begin{theorem}\label{the:CFM}
  Let $t > 0$ be a real number and let $A \colon [0,t] \to \GLpar{n}{\caln(G)} $ be a
  continuous piecewise smooth map.  Then
  \[ {\det}^{(2)}_{\caln(G)}(A(t)) = {\det}^{(2)}_{\caln(G)}(A(0)) \cdot
    \exp\left(\int_0^t \Re \left(\tr_{\caln(G)}\left(A(s)^{-1} \cdot
          \left.\frac{dA}{ds}\right|_s \right)\right)\,ds\right).
  \]
\end{theorem}




  \typeout{------ Section 4: Obstructions to fibring and virtually fibring -----------------}

  \section{Obstructions to fibring and virtually fibring coming from \texorpdfstring{$\phi$}{phi}-twisted
    \texorpdfstring{$L^2$}{L2}-torsion functions}%
\label{sec:Obstructions_to_fibring_and_virtually_fibring_coming_from_twisted_L2-torsion}


\subsection{The \texorpdfstring{$T_0$}{T0}-value and \texorpdfstring{$T_{\infty}$}{Too}-value and the property special}%
\label{subsec:T_0-value_T_infty-value_and_the_property_T-special}

Given a function $\rho \colon \IR_{>0} \to \IR$ and
$T_0,T_{\infty} \in \IR_{\ge 0} \amalg \{\infty\}$, we say that $\rho$ \emph{is
  $T_0$-special at $0$} if for some $d_0 \in \IZ$
\[
  \rho(t) = d_{0} \cdot \ln(t) \;\text{for}\; t < T_0
\]
holds and that $\rho$ \emph{is $T_{\infty}$-special at $\infty$} if for some
$d_{\infty} \in \IZ$
\[
  \rho(t) = d_{\infty} \cdot \ln(t) \;\text{for}\; T_{\infty} < t
\]
holds. Note that every $\rho$ is $0$-special at $0$ and $\infty$-special at $\infty$.

We say that $\rho$ is \emph{special} if it is both $T_0$-special at $0$ for some $T_0 > 0$
and $T_{\infty}$-special at infinity for some $T_{\infty} < \infty$.

  \begin{definition}[T-value]\label{def:T-values_and-special}
    Consider $\overline{\rho} \in \cale$ and $T_0,T_{\infty} \in \IR_{\ge 0} \amalg \{\infty\}$. We say that
    $\overline{\rho}$ is \emph{$T_0$-special} at $0$ if one (and hence every)
    representative $\rho$ of $\overline{\rho}$ is $T_0$-special and that $\overline{\rho}$
    is \emph{$T_{\infty}$-special} at $\infty$ if one (and hence every) representative
    $\rho$ of $\overline{\rho}$ is $T_{\infty}$-special at $\infty$. We call
    $\overline{\rho}$ special if one (and hence every) representative $\rho$ of
    $\overline{\rho}$ is special.

    We define the \emph{$T$-value} at $0$ and the \emph{$T$-value} at $\infty$ of
    $\overline{\rho}$ to be the elements in $\IR_{\ge 0} \amalg \{\infty\}$
    \begin{eqnarray*}
      T_0(\overline{\rho})
      & = &
            \sup\{T_0 \in \IR_{\ge 0} \amalg \{\infty\}\mid \overline{\rho}\; \text{is}\; T_0\text{-special at} \;0\};
      \\
      T_{\infty}(\overline{\rho})
      & = &
            \inf\{T_{\infty} \in \IR_{\ge 0} \amalg \{\infty\}\mid \overline{\rho}\; \text{is}\; T_{\infty}\text{-special at} \;\infty\}.
    \end{eqnarray*}
  \end{definition}

  Note  that $\overline{\rho} \in \cale$ is special if and only if both $T_0(\overline{\rho}) > 0$
  and $T_{\infty}(\overline{\rho}) < \infty$ hold.


\subsection{Specialness for mapping tori}%
\label{subsec:Chain_complex_version_for_mapping_tori}

Consider the following situation. Let $G$ be a group satisfying
Condition~\ref{basic_setup}~\ref{basic_setup:G}, i.e., $G$ is finitely generated and
residually finite and satisfies the Full Farrell--Jones Conjecture.  Let
$\phi \colon G \to \IZ$ be a surjective group homomorphism with kernel $K$.  Denote by
$i \colon K \to G$ the inclusion.  Consider an element $y \in G$ which is mapped under
$\phi \colon G \to \IZ$ to $1$.  Denote by $\gamma \colon K \to K$ the automorphism of $K$
sending $k$ to $yky^{-1}$. Let $D_*$ be a finite  free $\IZ[K]$-chain complex.
Denote by $\gamma^*D_*$ the finite free $\IZ[K]$-chain complex obtained from $D_*$
by restriction with $\gamma$.  Let $f_* \colon D_* \to \gamma^*D_*$ be a $\IZ[K]$-chain
map.  (We will not require that $f_*$ is a $\IZ[K]$-chain homotopy equivalence in
Theorem~\ref{the:mapping_tori_chain}.) We obtain a $\IZ[G]$-chain map
\[
  \widehat{f_*} \colon \IZ[G] \otimes_{\IZ[K]} D_* \to \IZ[G] \otimes_{\IZ[K]} D_*
\]
by sending $u \otimes v$ for $u \in \IZ[G]$ and $v \in D_n$ to
$u \otimes v - uy^{-1} \otimes f_n(v)$.  Recall from
Definition~\ref{def:chain_complex_version_of_a_mapping_torus} that the mapping torus
$T(f_*)_*$ of $f_*$ is defined to be $\cone(\widehat{f_*})$. Obviously it is a finite free
$\IZ[G]$-chain complex.

  \begin{theorem}\label{the:mapping_tori_chain}
    In the situation described above we get:

    \begin{enumerate}
    \item\label{the:mapping_tori_chain:det_cyclic} The mapping torus $T(f_*)_*$ is
      $L^2$-acyclic and hence its $\phi$-twisted reduced $L^2$-torsion function
      $\overline{\rho}^{(2)}(T(f_*)_*;\phi)$ is defined;

    \item\label{the:mapping_tori_chain:representative} Put
      $T = \max\{h(\widehat{f}_n) \mid n \in \IZ_{\ge 0}\}$.  Then there is a
      representative $\rho$ of $\overline{\rho}^{(2)}(T(f_*)_*;\phi)$ satisfying
      \[
        \rho(t) = 0 \quad \text{for}\; T < t.
      \]
    \end{enumerate}\smallskip
  \end{theorem}
  \begin{proof}~\ref{the:mapping_tori_chain:det_cyclic} This is the chain complex version
    of~\cite[Theorem~8.6 and Theorem~10.2]{Lueck(2018)}.
    \\[1mm]~\ref{the:mapping_tori_chain:representative} We use induction over the
    dimension of $D_*$. The induction beginning $\dim(D_*) = 0$ is done as follows. In
    this case $D_*$ is concentrated in dimension zero and given by a finitely generated
    free $\IZ[K]$-module $D_0$.  Choose a $\IZ[K]$-isomorphism
    $\mu \colon D_0 \xrightarrow{\cong} \IZ[K]^d$.  We obtain a $\IZ[K]$-homomorphism
    $f \colon \IZ[K]^d \to \gamma^*\IZ[K]^d$ by the composite
    $\IZ[K]^d \xrightarrow{\mu^{-1}} D_0 \xrightarrow{f_0} \gamma^*D_0
    \xrightarrow{\gamma^*\mu} \gamma^*\IZ[K]^d$.  Define the matrix
    $A \in \Mpar{d}{d}{\IZ[K]}$ by requiring that the image of the $i$-th element
    $(0,\ldots, 0, 1, 0, \ldots, 0)$ of the standard $\IZ[K]$-basis of $\IZ[K]^d$ is sent
    under $f$ to $(a_{i,1}, \ldots, a_{i,d})$. Consider the $1$-dimensional finitely generated  free
    $\IZ[G]$-chain complex $E_*$, whose first differential is the $\IZ[G]$-homomorphism
    $r_{I_d - y^{-1}A} \colon \IZ[G]^d \to \IZ[G]^d$ given by right multiplication with
    $I_d - y^{-1}A \in \Mpar{d}{d}{\IZ[G]}$.  We get a $\IZ[G]$-chain isomorphism
    $\widetilde{\mu}_* \colon T(f_*)_* \xrightarrow{\cong} E_*$ where
    $\widetilde{\mu}_i \colon T(f_*)_i = \IZ[G] \otimes_{\IZ[K]} D_0 \xrightarrow{\cong}
    \IZ[G]^d$ sends $u \otimes v$ to $u \cdot\mu(v)$.

    Consider $s \in \IR$ with $0 < s < \frac{1}{h(\widehat{f}_0)} = \frac{1}{h(y^{-1}A)}$.
    Then we obtain $0 \le h(r_{s\cdot y^{-1}A}) = h(s \cdot y^{-1}A) < 1$.  The operator
    given by the von Neumann series
    \[
      \sum_{i = 0}^{\infty} \big(r^{(2)} _{s \cdot y^{-1}A}\bigr)^i \colon L^2(G)^d \to L^2(G)^d
    \]
    is a well-defined $G$-equivariant isomorphism and is an inverse of
    \[r^{(2)}_{I_d - s\cdot y^{-1}A} \colon L^2(G)^d\to L^2(G)^d\] 
    by Lemma~\ref{lem:basic_properties_of_spectral_radius}~\ref{lem:basic_properties_of_spectral_radius:von_Neumann_series}.
    Fix $t \in \IR$ satisfying $h(\widehat{f}_0) < t$.  We get a smooth map from
    $[0,t^{-1}]$ to the space $ \GLpar{d}{\caln(G)}$ of invertible $G$-equivariant
    operators $L^2(G)^d \to L^2(G)^d$ equipped with the operator norm by sending
    $s \in [0,t^{-1} ]$ to $r^{(2)}_{I_d - s\cdot y^{-1}A}$.  Hence we conclude from
    Theorem~\ref{the:CFM} using the fact that $\tr_{\caln(G)}$ is continuous:
    \begin{eqnarray*}
      \lefteqn{{\det}_{\caln(G)}^{(2)}(r_{I_d - t^{-1}\cdot y^{-1}A})}
      & &
      \\
      & = & 
            {\det}_{\caln(G)}^{(2)}(r_{I_d - 0 \cdot y^{-1}A})
      \\
      & & \quad \quad \quad \quad \cdot \exp\left(\int_0^{t^{-1}}\Re\left(\tr_{\caln(G)}\left((r_{I_d - s\cdot y^{-1}A})^{-1}
            \circ  \frac{d}{ds}r_{I_d - sy^{-1}A}\right)\right) ds\right)
      \\
      &=&
          {\det}_{\caln(G)}^{(2)}(r_{I_d}) \cdot \exp\left(\int_0^{t^{-1}} \Re\left(\tr_{\caln(G)}
          \left(\left(\sum_{i=0}^\infty \bigl(r^{(2)}_{s \cdot y^{-1}A}\bigr)^i\right)
          \circ (-r^{(2)}_{y^{-1}A})\right)\right)ds\right)
      \\
      &=&
          \exp\left(\int_0^{t^{-1}}\Re\left(\tr_{\caln(G)}
          \left(\sum_{i=0}^\infty r^{(2)}_{s^i\cdot (y^{-1}A)^{i+1}}\right)\right)ds\right)
      \\
      &=&
          \exp\left(- \int_0^{t^{-1}}\Re\left(\sum_{i=0}^\infty (-s^i)
          \cdot \tr_{\caln(G)}\left( r^{(2)}_{(y^{-1}A)^{i+1}}\right)\right)ds\right).
    \end{eqnarray*}

    Each entry of the matrix $(y^{-1}A)^{i+1}$ is of the form
    $\sum_{j = 1}^k \lambda_j \cdot g_j$ for some $k \in \IZ_{\ge 0}$,
    $\lambda_j \in \IZ$, and $g_j \in G$ with $\phi(g_j) = -i - 1$.  Since
    $\phi(g_j) \neq 0$ implies $g_j \neq e_G$, we get
    $\tr_{\caln(G)}\left( r^{(2)}_{(y^{-1}A)^{i+1}}\right) = 0$ for $i \in \IZ_{\ge 0}$.
    Hence for $t$ with ${h(\widehat{f}_0)} < t $ we have
    \[ {\det}_{\caln(G)}^{(2)}(r^{(2)}_{I_d - t^{-1}\cdot y^{-1}A}) = 1.
    \]
    Hence
    \[\rho(t) = \ln\bigl({\det}_{\caln(G)}^{(2)}(r^{(2)}_{I_d - t^{-1}\cdot
      y^{-1}A})\bigr)\]
       is a representative for
    $\overline{\rho}(E_*;\phi) = \overline{\rho}(T(f_*)_*;\phi)$
    such that for
    $h(\widehat{f}_0) < t$ we have
    \[
      \rho(t) = 0.
    \]

    Next we deal with the induction step. Suppose that $D_*$ is $d$-dimensional.  Let
    $d[D_d]_*$ be the finite free $\IZ[K]$-chain complex concentrated in dimension $d$
    whose $d$th chain module is $D_d$. Let $D_*|_{d-1}$ be the finite free
    $(d-1)$-dimensional $\IZ[K]$-chain complex obtained from $D_*$ by restricting to
    dimension $(d-1)$.  Let $d[f_d]_* \colon d[D_d]_* \to \gamma^*d[D_d]_*$ and
    $f_*|_{d-1} \colon D_*|_{d-1} \to \gamma^*D_*|_{d-1}$ be the obvious $\IZ[K]$-chain
    maps induced by $f_*$. Then we obtain a commutative diagram of finite free
    $\IZ[K]$-chain complexes with exact rows
    \[\xymatrix@!C=32pt{0 \ar[r]
        & D_*|_{d-1} \ar[d]^{f_*|_{d-1}} \ar[r]^-{i_*} & D_* \ar[d]^{f_*} \ar[r]^-{p_*} &
        d[D_d]_* \ar[d]^{d[f_d]_*} \ar[r] & 0
        \\
        0 \ar[r] & \gamma^* D_*|_{d-1} \ar[r]^-{\gamma^*i_*} & \gamma^* D_*
        \ar[r]^-{\gamma^*p_*} & \gamma^* d[D_d]_* \ar[r] & 0.  }
    \]
    It induces a short exact sequence of finite free $\IZ[G]$-chain complexes
    \[
      0 \to T(f_*|_{d-1})_* \to T(f_*)_* \to T(d[f_d]_*)_* \to 0.
    \]
    By the induction hypothesis we can find representatives $\rho|_{d-1}$ of
    $\overline{\rho}( T(f_*|_{d-1})_*;\phi)$ and $\rho_d$ of
    $\overline{\rho}(T(d[f_d]_*)_*;\phi) = (-1)^d \cdot \overline{\rho}(T(0[f_d]_*)_*;\phi)$ satisfying
    \begin{eqnarray*}
      \rho|_{d-1}(t)  & =  & 0 \quad
                             \text{for}\;  \max\{h(\widehat{f}_n) \mid n \in \{0,1, \ldots, (d-1)\}\} < t;
      \\
      \rho_d(t)   & =  &  0 \quad \text{for}\;   h(\widehat{f_d}) < t.
    \end{eqnarray*}
    We conclude from Theorem~\ref{the:basic_properties_of_overline(rho)(C_ast;varphi)}~%
\ref{the:basic_properties_of_overline(rho)(C_ast;varphi):Additivity} that
    $\rho|_{d-1}(t) + \rho_d(t)$ is a representative of
    $\overline{\rho}^{(2)}(T(f_*)_*;\phi)$ and satisfies
    \[
      \rho|_{d-1}(t) + \rho_d(t)= 0 \quad \text{for} \max\{h(\widehat{f}_n) \mid n \in
      \{0,1, \ldots, d\} \} < t.
    \]
  \end{proof}

\begin{theorem}\label{the:finite_reduction_to_K}
  Let $C_*$ be a finite free $\IZ[G]$-chain complex. Denote by $i^*C_*$ the $\IZ[K]$-chain
  complex obtained from $C_*$ by restriction with $i$.  Suppose that there is a finite
  free $\IZ[K]$-chain complex $D_*$ which is $\IZ[K]$-chain homotopy equivalent to
  $i^*C_*$. Let $\chi_{\IZ[K]}(D_*) \in \IZ$ be the Euler characteristic of $D_*$. Then:

  \begin{enumerate}
  \item\label{the:finite_reduction_to_K:acyclic} The finite free $\IZ$[G]-chain complex
    $C_*$ is $L^2$-acyclic;
  \item\label{the:finite_reduction_to_K:representative} There are real numbers $T_0$ and
    $T_{\infty}$ with $0< T_0 \le T_{\infty}$ and a representative $\rho$ of
    $\overline{\rho}^{(2)}(C_*;\phi)$ satisfying
    \[
      \rho(t) =
      \begin{cases}
        0 & \text{for} \; t < T_0;
        \\
        \chi_{\IZ[K]}(D_*) \cdot \ln(t) & \text{for} \; T_{\infty} < t.
      \end{cases}
    \]
  \end{enumerate}
\end{theorem}
\begin{proof}
  We conclude from Lemma~\ref{lem:elementary_property_of_mapping_torus} that we are in the
  situation of Theorem~\ref{the:mapping_tori_chain} in the sense that there is a
  $\IZ[K]$-chain homotopy equivalence $f_* \colon D_* \to \gamma^*D_*$ such that $C_*$ and
  $T(f_*)_*$ are $\IZ[G]$-chain homotopy equivalent.  Hence by
  Theorem~\ref{the:mapping_tori_chain} there is $T_{\infty} > 0$ and a representative
  $\rho$ of $\overline{\rho}^{(2)}(T(f_*)_*;\phi)$ satisfying
  \[
    \rho(t) = 0 \quad \text{for}\; T_{\infty} < t.
  \]
  Now the same argument applies to $- \phi$ instead of $\phi$, since $\phi$ and $-\phi$
  have the same kernel.  So we find $T_0 > 0$ and a representative $\rho'$ of
  $\overline{\rho}^{(2)}(T(f_*)_*;-\phi)$ satisfying
  \[
    \rho'(t) = 0 \quad \text{for}\; T_0^{-1} < t.
  \]
  We conclude from Theorem~\ref{the:basic_properties_of_overline(rho)(C_ast;varphi)}~%
\ref{the:basic_properties_of_overline(rho)(C_ast;varphi):scaling} that there exists an
  integer $l$ such that for all $t \in \IR_{>0}$ we have
  $\rho'(t) + l \cdot \ln(t) = \rho(t^{-1}) $.  Hence we obtain
  \[
    \rho(t) =
    \begin{cases}
      -l \cdot \ln(t) & \text{for} \; t < T_0;
      \\
      0 & \text{for} \; T_{\infty} < t.
    \end{cases}
  \]
  If we replace $\rho(t)$ by $\rho(t) + l \cdot \ln(t)$, we get a representative $\rho$ of
  $\overline{\rho}^{(2)}(C_*;\phi)$ satisfying
  \[
    \rho(t) =
    \begin{cases}
      0 & \text{for} \; t < T_0;
      \\
      l \cdot \ln(t) & \text{for} \; T_{\infty} < t.
    \end{cases}
  \]
   
  It remains to show that $l = \chi_{\IZ[K]}(D_*)$ holds.
  By~\cite[Theorem~8.10]{Lueck(2018)} the limit
  \[\lim_{t \to \infty} \bigl(\rho(t) - \chi^{(2)}(L^2(G) \otimes_{\IZ[K]} D_*) \cdot
    \ln(t)\bigr) = \lim_{t \to \infty} \bigl((l - \chi^{(2)}(L^2(G) \otimes_{\IZ[K]} D_*))
    \cdot \ln(t)\bigr)
  \]
  exists. Hence we get $l = \chi^{(2)}(L^2(G) \otimes_{\IZ[K]} D_*) = \chi_{\IZ[K]}(D_*)$.
\end{proof}

Note that in Theorem~\ref{the:mapping_tori_chain} there exists in general no $T_0 > 0$ and
$C,D \in \IR$ such that $\rho(t) = C + D \cdot \ln(t)$ holds for $t < T_0$, unless $f_*$
is a chain homotopy equivalence.  Moreover, the induction step in the proof of
Theorem~\ref{the:mapping_tori_chain} does not work if we additionally assumed that
$f_*$ is a chain homotopy equivalence.


  \subsection{Property special is an obstruction to fibring over \texorpdfstring{$S^1$}{S1}}%
  \label{subsec:Property_special_is_an_obstruction_to_fibring_over_S_upper_1}

\begin{theorem}\label{the:Property_special_is_an_obstruction_to_fibring_and_virtually_fibring}
  Let $X$ be a connected finite $CW$-complex satisfying one of the following equivalent
  conditions:

  \begin{enumerate}

  \item There is a connected finite $CW$-complex $F$ such that $X$ is homotopy equivalent
    to the mapping torus $T_f$ for some selfhomotopy equivalence
    $f \colon F \xrightarrow{\simeq} F$;

  \item There exists a fibration $F \xrightarrow{i} E \xrightarrow{p} S^1$ such that $E$
    is homotopy equivalent to $X$ and $F$ is a connected finite $CW$-complex.

  \end{enumerate}

  Then:
  \begin{enumerate}
  
  \item\label{the:mapping_tori:L2-acyclic} The $L^2$-Betti number
    $b_n^{(2)}(\widetilde{X})$ vanishes for every $n \in \IZ_{\ge 0}$;

  \item\label{the:mapping_tori:overline(rho)} Let $\phi \colon \pi_1(X) \to \IZ$ be the
    surjective group homomorphism which is either associated to $T_f$
    or comes from $\pi_1(p)$ for the projection $p \colon E \to S^1$. Suppose that $\pi_1(X)$
    satisfies Condition~\ref{basic_setup}~\ref{basic_setup:G}.

    Then there are elements $T_0, T_{\infty} \in \IR_{> 0}$  and a
    representative $\rho \colon \IR_{> 0} \to \IR$ of
    $\overline{\rho}(\widetilde{X};\phi)$ satisfying
    \[
      \rho(t) = \begin{cases} 0 & \text{for} \; t < T_0;
        \\
        \chi(F) \cdot \ln(t) & \text{for} \; T_{\infty} < t.
      \end{cases}
    \]
    In particular $\overline{\rho}(\widetilde{X};\phi)$ is special.

  \end{enumerate}
\end{theorem}
\begin{proof}
  This is a direct consequence of Theorem~\ref{the:finite_reduction_to_K}.
\end{proof}

\begin{remark}\label{rem:T-values_and_Poincare_duality}If $M$ is a closed manifold, then
  $T_0(\overline{\rho}(\widetilde{M};\phi)) =
  T_{\infty}(\overline{\rho}(\widetilde{M};\phi))^{-1}$ follows from
  Theorem~\ref{the:Basic_properties_of_the_reduced_L2-torsion_function_for_universal_coverings}~%
\ref{the:Basic_properties_of_the_reduced_L2-torsion_function_for_universal_coverings:Poincare_duality}.
\end{remark}

\begin{remark}\label{rem:estimate_on_the_T-value}
  We get  explicit  bounds for the $T$-values of
  $\overline{\rho}(\widetilde{T_f};\phi)$ appearing in
  Theorem~\ref{the:Property_special_is_an_obstruction_to_fibring_and_virtually_fibring}
  in terms of the $L^1$-norms of the matrices associated to the
  $\IZ[\pi_1(f)]$-equivariant chain maps induced by a lift of $f$ or a homotopy inverse
  $f^{-1}$ to the universal covering $\widetilde{X}$ from
  Lemma~\ref{lem:basic_properties_of_spectral_radius}~%
\ref{lem:basic_properties_of_spectral_radius:inequality} and
  Theorem~\ref{the:mapping_tori_chain}.

  Let $f$ be a selfhomotopy equivalence of a closed connected orientable surface of genus
  $\ge 2$. Let $\phi \colon \pi_1(T_f) \to \IZ$ be the canonical group epimorphism.
  Denote by $e(\pi_1(f))$ the entropy of $\pi_1(f)$. Then we get
  $T_0(\overline{\rho}(\widetilde{T_f};\phi)) \le e(\pi_1(f))^{-1}$ and
  $e(\pi_1(f)) \le T_{\infty}(\overline{\rho}(\widetilde{T_f};\phi))$
  from~\cite[Theorem~8.5]{Dubois-Friedl-Lueck(2016)}.
\end{remark}

\begin{remark}\label{rem:twisted_phi_torsion_analytically}
  The $\phi$-twisted $L^2$-torsion should have an analytic counterpart. The hope is that
  one can compute it for hyperbolic manifolds since the heat kernel on the universal
  coverings of hyperbolic manifolds is well-understood and one may be able to extend this
  to the twisted setting, so that one can decide whether the property special is
  satisfied or not.

  This strategy would have the advantage that one can really exploit hyperbolicity which
  is rather difficult in the more algebraic approaches so far.
\end{remark}

\begin{remark}\label{exa:graph_manifold}
  It is explained in~\cite[Appendix~A]{Hughes-Leary-Lueck(2026not_fibre)} that the only
  class of closed $3$-manifolds, where we cannot decide in general whether they virtually
  fibre over $S^1$ is the class of a graph manifold $M$ whose fundamental group is not
  {\RFRS}.

  We get no new obstruction to virtual fibring over $S^1$ from the invariants defined in
  this paper for a graph manifold $M$.  Namely, we conclude
  from~\cite[Theorem~1.2]{Dubois-Friedl-Lueck(2016)} that
  $\overline{\rho}(\widetilde{N};\phi)$ is special for every finite covering $N$ of $M$
  and every $\phi \in H^1(N;\IZ)$, since every finite covering of a graph manifold is
  again a graph manifold.
\end{remark}


\subsection{Application to groups}\label{sec:Application_to_groups}

Note that there is an obvious version of
Theorem~\ref{the:Property_special_is_an_obstruction_to_fibring_and_virtually_fibring}
for groups.  Namely, let $G$ be a group such that $G$ satisfies
Condition~\ref{basic_setup}~\ref{basic_setup:G} and is of type FP, i.e., there exists a
finite projective $\IZ[G]$-resolution $P_*$ of the $\IZ[G]$-module $\IZ$ which is given by
the abelian group $\IZ$ with trivial $G$-action.  Note that this implies that $G$ is
torsionfree.  Let $\phi \colon G \to \IZ$ be a  group homomorphism.

Since $G$ satisfies the Full Farrell--Jones Conjecture, we can choose $P_*$ to be a finite
free $\IZ[G]$-chain complex and can define
\begin{equation}
  \overline{\rho}^{(2)}(G;\phi) := \overline{\rho}^{(2)}(P_*;\phi) \in \cale.
  \label{rho_upper_(2)(G;varphi)}
\end{equation}
This is independent of the choice of $P_*$ because of
Theorem~\ref{the:basic_properties_of_overline(rho)(C_ast;varphi)}~%
\ref{the:basic_properties_of_overline(rho)(C_ast;varphi)_homotopy_invariance}.

\begin{lemma}\label{lem:T-positive_as_obstruction_for_fibring_for_groups}
  Suppose that the kernel $K$ of  the non-trivial $\phi \colon G \to \IZ $ is of type FP$_{\infty}$, i.e., there exists a
  (not necessarily finite-dimensional) finitely generated projective $\IZ[K]$-resolution
  of the $\IZ[K]$-module $\IZ$.

  Then $\overline{\rho}^{(2)}(G;\phi)$ is special.
\end{lemma}
\begin{proof}
  Since there exists both a finite-dimensional projective $\IZ[K]$-resolution and a
  finitely generated projective $\IZ[K]$-resolution of the $\IZ[K]$-module $\IZ$, there
  exists a finite projective $\IZ[K]$-resolution $P_*$ of the $\IZ[K]$-module $\IZ$, see for
  instance~\cite[Proposition~11.19 on page~271]{Lueck(1989)}. Since $G$ and hence $K$
  satisfies the Full Farrell--Jones Conjecture, we can arrange $P_*$ to be  finite free.
  Theorem~\ref{the:finite_reduction_to_K} implies that   $\overline{\rho}^{(2)}(G;\phi)$ is special.
\end{proof}


\typeout{--------- Section 9: The special case of an infinite cyclic group  -----------------}

\section{The special case of an infinite cyclic group}%
\label{sec:The_special_case_of_an_infinite_cyclic_group}

In this section we consider as an illustration the special case $G = \IZ$ and
$\phi = \id_{\IZ}$.


\subsection{Fuglede--Kadison determinant over \texorpdfstring{$\IZ$}{Z}}%
\label{subsec:Fuglede--Kadison_determinants_over_IZ}

Consider an element $p(z) \in \IZ[\IZ] = \IZ[z,z^{-1}]$ with $p \not = 0$. It can be
written uniquely in the form
\[p(z) = \sum_{n = n_-}^{n_+} s_n \cdot z^n
\]
for $n_-,n_+,s_n \in \IZ$  with $n_- \le n_+$, $s_{n_-} \not = 0$, and
$s_{n_+} \not=0$.  The \emph{lower degree $\deg_-(p)\in \IZ$} is defined to be $n_-$,
whereas the \emph{upper degree $\deg_+(p)\in \IZ$} is defined to be $n_+$. Define the
\emph{width} $\width(p) \in \IZ_{\ge 0}$ of $p$ to be $\deg_+(p) - \deg_-(p)$. The
\emph{lower leading coefficient} $\lead_-(p) \in \IZ$ is defined to be $s_{n_-}$, whereas
the \emph{upper leading coefficient} $\lead_+(p) \in \IZ$ is defined to be $s_{n_+}$.  We
assign to $p$ an element
\begin{equation}
  \overline{\rho}(p) \in \cale
  \label{overline(rho)(p)}
\end{equation}
as follows.  Consider the $1$-dimensional $\IZ[\IZ]$-chain complex $C_*$ whose first
differential is the $\IZ[\IZ]$-map $r_p \colon \IZ[\IZ] \to \IZ[\IZ]$ given by
multiplication with $p$. Then $C_*$ is $L^2$-acyclic by~\cite[Lemma~1.34 on
page~35]{Lueck(2002)} and we define $\overline{\rho}(p)$ to be the $\id_{\IZ}$-twisted
$L^2$-torsion function $\overline{\rho}(C_*;\id_{\IZ})$ of $C_*$,
see~\eqref{reduced_twisted_L_upper_torsion_function_chain}.

We can write $p$ in $\IC[\IZ] = \IC[z,z^{-1}]$ with $p \not= 0$ uniquely to be
\[
  p(z) = C \cdot z^k \cdot \prod_{i = 1}^l (z- a_i)
\]
for $k \in \IZ$, $C \in \IZ\setminus \{0\}$, $l \in \IZ_{\ge 0}$, and $a_i \in \IC \setminus \{0\}$. We conclude
\[
  \ln\big( {\det}_{\caln(\IZ)}(r_p^{(2)} \colon L^2(\IZ) \to L^2(\IZ))\bigr) = \ln(|C|) +
  \sum_{\substack{i \in \{1,2, \ldots , l\}\\ |a_i| \ge 1}} \ln(|a_i|)
\]
from~\cite[(3.23) on page~136]{Lueck(2002)}. Note that the right hand side of the equation
above is just the logarithm of the Mahler measure of $p$.

For $t \in \IR_{> 0}$ define $p_t \in \IC[z,z^{-1}]$ by $p_t(z) = p(tz)$. We can write
\[
  p_t(z) = C \cdot t^{k+l} \cdot z^k \cdot \prod_{i = 1}^l \left( z- \frac{a_i}{t}\right).
\]
Then the following function $\rho \colon \IR_{> 0} \to \IR$ is a representative for
$\overline{\rho}(p) \in \cale$
\begin{eqnarray*}
  \rho(t)
  & = &
        \ln\bigl({\det}_{\caln(\IZ)}(r_{p_t}^{(2)} \colon L^2(\IZ) \to L^2(\IZ))\bigr)
  \\
  & = &
        \ln(|C|) + (k+l) \cdot \ln(t) + \sum_{\substack{i \in \{1,2, \ldots , l\}\\ |a_i/t|\ge 1}} \ln(|a_i/t|)
  \\
  & = &
        \ln(|C|) + \sum_{\substack{i \in \{1,2, \ldots, l\}\\ |a_i| \ge  t}}
  \ln(|a_i|) + (k + |\{i \in \{1,2, \ldots, l\}\mid  |a_i| < t \}|) \cdot \ln(t).
\end{eqnarray*}
This implies
\begin{equation}
  \rho(t) = \begin{cases}
    \ln(|C|) +  \sum_{i = 1}^l \ln(|a_i|) + k \cdot \ln(t) & t < \min\{|a_i| \mid i = 1,2, \ldots, l\};
    \\
    \ln(|C|) + (k+l) \cdot \ln(t) & t > \max\{|a_i| \mid i = 1,2, \ldots, l\}.
    \label{asymptotics_for_rho_p(t)}
  \end{cases}
\end{equation}
Write $\zero(p) = k$ for the \emph{multiplicity} of $0$ as a zero of $p$.  If we put
\[
  q(z) = \prod_{i = 1}^l \left( z- a_i\right),
\]
then $q \in \IQ[z]$ and is uniquely determined by
\[
  p(z) = \lead_+(p) \cdot z^{\zero(p)} \cdot q(z).
\]
We have
\[\lead_-(p) = \lead_+(p) \cdot q(0) = \lead_+(p) \cdot \prod_{i = 1}^l (-a_i).
\]
Then equation~\eqref{asymptotics_for_rho_p(t)} becomes
\begin{equation}
  \rho(t) = \begin{cases}
    \ln(|\lead_-(p)|) + \zero(p)  \cdot \ln(t) & t < \min\{|a_i| \mid i = 1,\ldots, l\};
    \\
    \ln(|\lead_+(p)|) + (\zero(p)+\width(p)) \cdot \ln(t) & t > \max\{|a_i| \mid i = 1,\ldots, l\}.
    \label{asymptotics_for_rho_p(t)_rewritten}
  \end{cases}
\end{equation}

For the involution $\bar\ \colon \IZ[\IZ] \to \IZ[\IZ]$ defined by
$\overline{\sum_{n \in \IZ} s_n \cdot z^n} = \sum_{n \in \IZ} s_n \cdot z^{-n}$ we get
\begin{eqnarray*}
  \deg_-(\overline{p})
  &=&
      -\deg_+(p);
  \\
  \deg_+(\overline{p})
  &=&
      -\deg_-(p);
  \\
  \width(\overline{p}) 
  & = & 
        \width(p);
  \\
  \lead_+(\overline{p})
  & = & 
        \lead_-(p);
  \\
  \lead_-(\overline{p})
  & = & 
        \lead_+(p);
  \\
  \zero(\overline{p})
  & = &
        -\zero(p) - \width(p).
\end{eqnarray*}
Since we have $\overline{p}_t = \overline{p_{t^{-1}}}$ and
$\det_{\caln(\IZ)}(r_{\overline{p_{t^{-1}}}}) = \det_{\caln(\IZ)}(r_{p_{t^{-1}}})$, the
function $t \mapsto \rho(t^{-1})$ is a representative of $\overline{\rho}(\overline{p})$.

\begin{remark}
	The $\id_\IZ$-twisted $L^2$-torsion functions of free finite $\IZ[\IZ]$-chain complexes at 0 and $\infty$ take the form $C+D\cdot\ln(t)$.  Christopher Smyth, in personal communication, showed us that this fails already for $\IZ^2$, the epimorphism $\phi\colon \IZ^2 \to \IZ$ sending $z_1$and $z_2$ to $1$, and the operator $r_q$ corresponding to $q = 1 + z_1 + z_2 \in \IZ[\IZ^2]$, where $z_1 = (1,0)$ and $z_2 = (0,1)$. Smyth's argument used a computation of Mahler measures due to Maillot and  Cassaigne \cite[Proposition~7.3.1]{Maillot(2000)}.  We do not reproduce the argument here.
\end{remark}


  \subsection{Alexander polynomials}\label{subsec:Alexander_polynomials}

  Let $S$ be a \UFD, i.e., a commutative Noetherian unique factorization domain. Examples
  are principal ideal domains such as $\IZ$ or any field. Note that the group ring $S[A]$ is a
  \UFD\ if $S$ is a \UFD\ and $A$ a finitely generated free abelian group. Let $M$ be a
  finitely generated $S$-module. Then we can find an $(a,b)$-matrix $C \in \Mpar{a}{b}{S}$
  for $a,b > 0$ such that the cokernel of the $S$-homomorphism $r_C \colon S^a \to S^b $
  given by right multiplication with $C$ is $S$-isomorphic to $M$.  Let $E(M)$ be the
  \emph{elementary ideal of $C$}, i.e., the ideal in $S$ generated by the $(b,b)$-minors
  of $C$, if $b \le a$ and  $\{0\}$ if $b > a$. It turns out that $E(M)$ depends only
  on $M$ and not on the choice of $C$.  Since $S$ is a \UFD, there exists a unique
  smallest principal ideal $I$ of $S$ which contains $E(M)$. A generator of $I$ is called
  the \emph{order} $\ord(M) \in S$ of $M$.  The order is only well-defined up to
  multiplication with units in $S$. A representative of $\ord(M)$ is given by the greatest
  common divisor of the $(b,b)$-minors of $C$, if $b \le a$. Note that using Cramer's rule
  one can show that each of the $(b,b)$-minors and hence $\ord(M)$ annihilates $M$.

  Let $S^{(0)}$ be the quotient field of the commutative integral domain $S$.  The
  elementary proof of the next result is left to the reader.

\begin{lemma}\label{lem:basics_properties_of_S-modules_UFD}
  The following statements are equivalent for a finitely generated $S$-module $M$:
  \begin{enumerate}
  \item $M$ is an $S$-torsion module;
  \item $S^{(0)} \otimes_S M = \{0\}$;
  \item $\dim_{S^{(0)}}(S^{(0)} \otimes_S M) = 0$;
  \item $\ord(M) \not= 0$.
  \end{enumerate}
\end{lemma}

We will be interested in the case $S = \IZ[\IZ] = \IZ[s,s^{-1}]$.  Note that the
multiplicative group of units in $\IZ[\IZ]$ is $U := \{\pm s^m \mid m \in \IZ\}$. Then for
any finitely generated $\IZ[\IZ]$-module $M$ we have its order
\[
  \ord(M) \in \IZ[\IZ]/U.
\]
\begin{lemma}\label{lem:computing_the_order}
  Let $M$ be a $\IZ[\IZ]$-module such that $M$ is finitely generated as abelian group.

  Then  $M/\tors(M)$ is finitely generated free as abelian group.  In
  particular we can consider the characteristic polynomial
  $\det(s - l_g) \in \IZ[s] \subseteq \IZ[s,s^{-1}] = \IZ[\IZ]$ of the automorphism
  $l_g \colon M/\tors(M) \to M/\tors(M) $ given by multiplication with the standard generator  $g$ of
  $\IZ$.

  Then we get in $\IZ[\IZ]/U$
  \[\ord(M) = \ord(M/\tors(M)) = \det(s - l_g).
  \]
\end{lemma}
\begin{proof} We get $\ord(M) = \ord(M/\tors(M))$
  from~\cite[page~2]{Silver-Williams(2012)}.

  We have the exact sequence of $\IZ[\IZ]$-modules 
  \[0 \to \IZ[\IZ] \xrightarrow{l_{g-1}} \IZ[\IZ] \xrightarrow{\epsilon} \IZ \to 0,
  \]
  where $\IZ$ is equipped with the trivial $\IZ$-action, $l_{g-1}$ is multiplication with
  $g-1$, and $\epsilon$ is the augmentation homomorphism sending $g$ to $1$.  Given a
  $\IZ[\IZ]$-module $N$ whose underlying abelian group is finitely generated free, we
  obtain the exact sequence of $\IZ[\IZ]$-modules
  \[0 \to \IZ[\IZ] \otimes_{\IZ} N \xrightarrow{l_{g-1}\otimes_{\IZ} \id_N} \IZ[\IZ]
    \otimes_{\IZ} N \xrightarrow{\epsilon \otimes_{\IZ} \id_N} \IZ \otimes_{\IZ} N = N \to
    0,
  \]
  where $\IZ$ acts on the tensor products diagonally and $\IZ[\IZ] \otimes_{\IZ} N$ is
  finitely generated free as $\IZ[\IZ]$-module.  Now $\ord(N) = \det(s - l_g)$ follows
  from the definitions.
\end{proof}

Let $\IZ[\IZ]^{(0)} = \IZ[z,z^{-1}]^{(0)}$ be the quotient field of the commutative
integral domain $\IZ[\IZ] = \IZ[z,z^{-1}]$.  Let $C_*$ be a finite free $\IZ[\IZ]$-chain
complex. Then the $\IZ[\IZ]$-homology groups $H_i(C_*)$ are finitely generated and are
non-zero only for finitely many values of $i$.

\begin{lemma}\label{lem:L_upper_2-acyclic_and_Alexander}
  The following assertions are equivalent for a finite free $\IZ[\IZ]$-chain complex
  $C_*$.

  \begin{enumerate}
  \item\label{lem:L_upper_2-acyclic_and_Alexander:L_upper_acyclic} $C_*$ is $L^2$-acyclic,
    i.e., all $L^2$-Betti numbers $b_n^{(2)}(C_*;\caln(\IZ))$ vanish;
  \item\label{lem:L_upper_2-acyclic_and_Alexander:dimension} We have
    $\dim_{\IZ[\IZ]^{(0)}}(\IZ[\IZ]^{(0)} \otimes_{\IZ[\IZ]}H_n(C_*)) = 0$ for every
    $n \in \IZ_{\ge 0}$;

  \item\label{lem:L_upper_2-acyclic_and_Aleander:order} We have $\ord(H_n(C_*)) \not = 0$
    for every $n \in \IZ_{\ge 0}$.
  \end{enumerate}
\end{lemma}
\begin{proof}\ref{lem:L_upper_2-acyclic_and_Alexander:L_upper_acyclic}
  $\Longleftrightarrow$~%
\ref{lem:L_upper_2-acyclic_and_Alexander:dimension} See~\cite[Lemma~1.34 on
  page~35]{Lueck(2002)}.
  \\[1mm]\ref{lem:L_upper_2-acyclic_and_Alexander:dimension}$\Longleftrightarrow$~%
\ref{lem:L_upper_2-acyclic_and_Aleander:order} This follows from
  Lemma~\ref{lem:basics_properties_of_S-modules_UFD}.
\end{proof}

Hence we can define for a finite free $\IZ[\IZ]$-chain complex $C_*$ which is
$L^2$-acyclic elements
\begin{eqnarray}
  \Delta_i(C_*)
  & = &
        \ord(H_i(C_*)) \quad \in \IZ[\IZ]/U;
        \label{Delta_i(C_ast)}
  \\
  \rho_{\Delta}(C_*)
  & = &
        \prod_{i \ge 0} \Delta_i(C_*)^{(-1)^i} \in \IZ[\IZ]^{(0)}/U,
        \label{rho_Delta(C_ast)}
\end{eqnarray}
where we call $\Delta_i(C_*)$ the \emph{i-th Alexander polynomial} of $C_*$.

One can also define the universal $L^2$-torsion for a finite free $\IZ[\IZ]$-chain complex
$C_*$ which is $L^2$-acyclic, see~\cite[Section~2 and
Example~3.3]{Friedl-Lueck(2017universal)}
\begin{equation}
  \rho^{(2)}_u(C_*) \in  \Wh^w(\IZ)
  \label{rho_upper(C_ast)}
\end{equation}

\begin{lemma}\label{lem:Alexander_polynomial_and_universal_L_upper_2-torsion}
  Let $C_*$ be a finite free $\IZ[\IZ]$-chain complex which is $L^2$-acyclic. Then
  there is an isomorphism
  \[
    \omega \colon \Wh^w(\IZ) \xrightarrow{\cong} \IZ[s,s^{-1}]^{(0)}/U
  \]
  which sends $\rho^u(C_*)$ to $\rho_{\Delta}(C_*)$.
\end{lemma}
\begin{proof} This follows from~\cite[(3.4)]{Friedl-Lueck(2017universal)}
  and~\cite[Theorem~2.6 on page~11 and Theorem~4.7 on page~20]{Turaev(2001)}.
\end{proof}
Lemma~\ref{lem:Alexander_polynomial_and_universal_L_upper_2-torsion} together
with~\cite[Section~3.4]{Friedl-Lueck(2017universal)} implies

\begin{theorem}\label{them:Alexander_polynomial_and_twisted_L_upper_2-torsion}
  Let $C_*$ be a finite free $\IZ[\IZ]$-chain complex which is $L^2$-acyclic. Then we get
  for the $\id_{\IZ}$-twisted $L^2$-torsion function $\overline{\rho}(C_*;\id_{\IZ} )$ of
  $C_*$
  \[
    \overline{\rho}(C_*;\id_{\IZ} ) = \sum_{i \ge 0} (-1)^i \cdot
    \overline{\rho}(\ord(H_i(C_*));\id_{\IZ}),
  \]
  where for $\overline{x} \in \IZ[\IZ]/U$ we define
  $\overline{\rho}(\overline{x};\id_{\IZ})$ to be the $\id_{\IZ}$-twisted $L^2$-torsion
  function $\overline{\rho}(\el[1,x]_*;\id_{\IZ} )$ of the $1$-dimensional free
  $\IZ[\IZ]$-chain complex $\el[1,x]_*$ whose first differential is given by
  $\IZ[\IZ] \to \IZ[\IZ], \;a \mapsto ax$.
\end{theorem}

Recall that we have computed $\overline{\rho}(\ord(H_i(C_*)))$ explicitly in
Subsection~\ref{subsec:Fuglede--Kadison_determinants_over_IZ}. In particular there exists
$T \in \IR_{> 0}$ such that we get a representative for $\overline{\rho}(C_*;\id_{\IZ})$
by
\begin{equation}
  \rho(t) =
  \sum_{i \ge 0} (-1)^i \cdot \ln(|\lead_-(p_i)|)
  + \biggl(\sum_{i \ge 0} (-1)^i \cdot  \zero(p_i) \biggr) \cdot \ln(t)
  \label{asymptotics_for_rho_p(t)_Alexander_small_t}
\end{equation}
for $t < T$ and
\begin{equation}
  \rho(t) =
  \sum_{i \ge 0}  (-1)^i \cdot  \ln(|\lead_+(p_i)|) + \biggl(\sum_{i \ge 0} (-1)^i \cdot \bigl(\zero(p_i)
    + \width(p_i)\bigr) \biggr)\cdot \ln(t)
  \label{asymptotics_for_rho_p(t)_Alexander_large_t}
\end{equation}
for $t > T^{-1}$, where $p_i$ is a representative for $\Delta_i(C_*) \in \IZ[\IZ]/U$ and
$q_i \in \IQ[z]$ is determined by the properties $q_i(0) \not= 0$ and
$p_i(z) = \lead_+(p_i) \cdot z^{\zero(p_i)} \cdot q_i(z)$.


\subsection{Leading coefficients and mapping tori}%
\label{subsec:leading_coefficients_and_mapping_tori}

  \begin{lemma}\label{lem:sufficient_condition_for_fibring_in_case_Z}
    The following assertions are equivalent for a finite free $\IZ[\IZ]$-chain complex
    $C_*$.

    \begin{enumerate}
    \item\label{lem:sufficient_condition_for_fibring_in_case_Z:mapping_torus} The
      $\IZ[\IZ]$-chain complex $C_*$ is $\IZ[\IZ]$-chain homotopy equivalent to the
      mapping torus $T(f_*)_*$ of a $\IZ$-chain map $f_* \colon D_* \to D_*$ of a finite
      free $\IZ$-chain complex $D_*$;

    \item\label{lem:sufficient_condition_for_fibring_in_case_Z:finite_free_chain} The
      $\IZ$-chain complex $i^*C_*$, which is obtained from $C_*$ by restriction with the
      inclusion $i \colon \IZ \to \IZ[\IZ]$, is $\IZ$-chain homotopy equivalent to a
      finite free $\IZ$-chain complex $D_*$;

    \item\label{lem:sufficient_condition_for_fibring_in_case_Z:homology_fin_gen} The
      homology groups $H_n(i^*C_*)$ are finitely generated as abelian group for $n \ge 0$
      and non-trivial only for finitely many $n \ge 0$;

    \item\label{lem:sufficient_condition_for_fibring_in_case_Z:leading_term} For every
      $n \ge 0$ the element $\Delta_n(C_*) \in \IZ[\IZ]/U$ is non-zero and both leading
      coefficients $\lead_-(\Delta_n(C_*))$ and $\lead_+(\Delta_n(C_*))$ belong to $\{\pm 1\}$;

    \item\label{lem:sufficient_condition_for_fibring_in_case_Z:Novikov} The homology
      groups $H_n(\IZ((t)) \otimes_{\IZ[t,t^{-1}]} C_*)$ and
      $H_n(\IZ((t^{-1})) \otimes_{\IZ[t,t^{-1}]} C_*)$ vanish for all $n \ge 0$.
    \end{enumerate}
  \end{lemma}
  \begin{proof}~\ref{lem:sufficient_condition_for_fibring_in_case_Z:mapping_torus}
    $\Longleftrightarrow$~\ref{lem:sufficient_condition_for_fibring_in_case_Z:finite_free_chain}
    This follows from Lemma~\ref{lem:elementary_property_of_mapping_torus}.
    \\[1mm]~\ref{lem:sufficient_condition_for_fibring_in_case_Z:finite_free_chain}
    $\Longleftrightarrow$~\ref{lem:sufficient_condition_for_fibring_in_case_Z:homology_fin_gen}
    Since $\IZ$ is Noetherian, we
    get~\ref{lem:sufficient_condition_for_fibring_in_case_Z:finite_free_chain}
    $\implies$~\ref{lem:sufficient_condition_for_fibring_in_case_Z:homology_fin_gen}.
    The other
    implication~\ref{lem:sufficient_condition_for_fibring_in_case_Z:homology_fin_gen}
    $\implies$~\ref{lem:sufficient_condition_for_fibring_in_case_Z:finite_free_chain} is
    proved as follows.  Let $N$ be the dimension of $C_*$. Since $\IZ$ is a principal
    ideal domain, we can choose for any $n \ge 0$ a finite free $1$-dimensional chain
    complex $E[n]_*$ such that $H_1(E[n]_*)$ vanishes and
    $H_n(i^*C_*) \cong_{\IZ} H_0(E[n]_*)$ holds.  Then one easily constructs a $\IZ$-chain
    map $h_* \colon \bigoplus_{n = 0}^N \Sigma^nE[n]_* \to i^*C_*$ which induces
    isomorphisms on all homology groups. Since its source and target are projective, $h_*$
    is a $\IZ$-chain homotopy equivalence.  Obviously the source is a finite free
    $\IZ$-chain complex.
    \\[1mm]~\ref{lem:sufficient_condition_for_fibring_in_case_Z:homology_fin_gen}
    $\implies$~\ref{lem:sufficient_condition_for_fibring_in_case_Z:leading_term}
    This follows from Lemma~\ref{lem:computing_the_order}.
    \\[1mm]~\ref{lem:sufficient_condition_for_fibring_in_case_Z:leading_term}
    $\implies$~\ref{lem:sufficient_condition_for_fibring_in_case_Z:homology_fin_gen}
    This follows from~\cite[Corollary~2.5]{Silver-Williams(2012)}. 
    \\[1mm]~\ref{lem:sufficient_condition_for_fibring_in_case_Z:finite_free_chain}
    $\Longleftrightarrow$~\ref{lem:sufficient_condition_for_fibring_in_case_Z:Novikov}
    This follows from~\cite[Theorem~2]{Ranicki(1995Novikovring)}.
  \end{proof}

\begin{lemma}\label{lem:leading_cowfficient_mapping_tori_over_Z}
  Let $C_*$ be a finite free $\IZ[\IZ]$-chain complex. Suppose that $C_*$ is
  $\IZ[\IZ]$-chain homotopy equivalent to the mapping torus $T(f_*)_*$ of a $\IZ$-chain
  map $f_* \colon D_* \to D_*$ of a finite free $\IZ$-chain complex $D_*$.

  Then $\ord(H_i(C_*)) \not = 0$ holds for every $i \in \IZ_{\ge 0}$ and we get
  \[
    \prod_{i \ge 0} |\lead_+(\Delta_i(C_*))|^{(-1)^{i}} = 1.
  \]
\end{lemma}
\begin{proof}
  We get from Theorem~\ref{the:mapping_tori_chain},
  Lemma~\ref{lem:Alexander_polynomial_and_universal_L_upper_2-torsion},
  and\eqref{asymptotics_for_rho_p(t)_Alexander_large_t} that $\ord(H_i(C_*)) \not = 0$
  holds for every $i \in \IZ_{\ge 0}$ and that there are $T_{\infty} \in \IR_{> 0}$ and
  $m\in \IZ$ satisfying
  \[
    \sum_{i \ge 0} (-1)^i \cdot \ln\bigl(|\lead_+(\Delta_i(C_*))|\bigr) = m \cdot \ln(t)
  \]
  for $t > T_{\infty}$. This implies
  $\sum_{i \ge 0} (-1)^i \cdot \ln\bigl(|\lead_+(\Delta_i(C_*))|\bigr) = 0$.
\end{proof}

The conditions appearing in Lemma~\ref{lem:leading_cowfficient_mapping_tori_over_Z} are
necessary but not sufficient for fibring as the next example illustrates.

\begin{example}\label{exa:not_sufficient}
  Let $C_*$ be the $2$-dimensional finite free $\IZ[\IZ]$-chain complex given in
  dimensions $2,1,0$ by
  \[\IZ[\IZ] \xrightarrow{\begin{pmatrix} 2 \\ 0 \end{pmatrix}} \IZ[\IZ] \oplus\IZ[\IZ]
    \xrightarrow{\begin{pmatrix} 0 & 2 \end{pmatrix}} \IZ[\IZ]
  \]
  
  Then the $\IZ[\IZ]$-module $H_i(C_*)$ is $\IZ[\IZ]$-isomorphic to the cokernel of
  $2 \cdot \id \colon \IZ[\IZ] \to \IZ[\IZ]$ and hence $\IZ/2[\IZ]$ for $i = 0,1$ and
  vanishes otherwise. Hence a representative for $\Delta_i(C_*)$ is the constant
  polynomial $2$.  We conclude that $\ord(H_i(C_*)) \not = 0$ holds for every
  $i \in \IZ_{\ge 0}$ and we have $\prod_{i \ge 0} |\lead_+(\Delta_i(C_*))|^{(-1)^{i}} = 1$.

  On the other hand $C_*$ considered as $\IZ$-chain complex is not $\IZ$-chain homotopy
  equivalent to a finite free $\IZ$-chain complex by
  Lemma~\ref{lem:sufficient_condition_for_fibring_in_case_Z} as $H_i(C_*)$ is not
  finitely generated as a $\IZ$-module.
\end{example}

Theorem~\ref{them:Alexander_polynomial_and_twisted_L_upper_2-torsion} and
Lemma~\ref{lem:leading_cowfficient_mapping_tori_over_Z} illustrate that the $\phi$-twisted
$L^2$-torsion indeed yields an obstruction for fibring and generalizes the obstruction given
by the alternating product of leading coefficients of the Alexander polynomials from infinite cyclic coverings to
arbitrary coverings.

\subsection{Leading coefficients and the matrix 3g+2e}
\label{subsec:Back_to_the_non-aspherical_counterexample}

As an illustration we go back to the example appearing
in~\cite[Section~3]{Hughes-Leary-Lueck(2026not_fibre)}.

Consider a group $G$ and an element $g \in G$ such that for some subgroup $G'$ of finite
index we have $g \in G'$ and its image under the canonical projection $G' \to H_1(G')_f$
is non-trivial.  Let $C_*$ be the $1$-dimensional finite free $\IZ[G]$-chain complex whose
first differential $c_1 \colon \IZ[G] \to \IZ[G]$ is given by right multiplication with
$3g+2e$. Then $C_*$ is $L^2$-acyclic.  We want to show for any group epimorphism
$\phi \colon G \to \IZ$ with kernel $K$ that the $\IZ[K]$-chain complex $i^*C_*$ is not
$\IZ[K]$-chain homotopy equivalent to a finite projective $\IZ[K]$-chain complex for
$i \colon K \to G$ the inclusion.

Suppose the contrary.  Let $j \colon \{1\} \to \IZ$ be the inclusion. Let
$\pr \colon K \to \{1\}$ be the projection.  Then there is an isomorphism of
$\IZ$-chain complexes
\[
  \pr_* i^* C_* = \IZ \otimes_{\IZ[K]} i^*C_* \xrightarrow{\cong} j^*\phi_*C_* = \IZ[\IZ]
  \otimes_{\IZ[G]} C_*
\]
sending $1 \otimes x$ to $z^0 \otimes x$, where $z$ is the standard generator of $\IZ$.
If $u \in G$ is any preimage of $z$ under $\phi$, then the inverse sends $z^m \otimes x$
to $1 \otimes u^mx$. Let $D_*$ be the $1$-dimensional free $\IZ[\IZ]$-chain complex whose
first differential is multiplication with $3z^{\phi(g)} + 2$ for $z \in \IZ$ the standard generator.
Since $D_* = \phi_*C_*$, we
conclude that $j^*D_*$ is $\IZ$-chain homotopy equivalent to a finite projective
$\IZ$-chain complex $E_*$.  As $\IZ$ is a principal ideal domain, $E_*$ is a finite
free $\IZ$-chain complex.  Hence there is an endomorphism $f_* \colon E_* \to E_*$
such that $D_*$ is $\IZ[\IZ]$-chain homotopy equivalent to the mapping torus $T(f_*)_*$,
see~\cite[Lemma~3.3]{Hughes-Leary-Lueck(2026not_fibre)}.
Lemma~\ref{lem:elementary_property_of_mapping_torus} and 
Lemma~\ref{lem:leading_cowfficient_mapping_tori_over_Z} imply
$\lead_+(\Delta_1(D_*))^{-1} \cdot \lead_+(\Delta_0(D_*)) = 1$.  Since $H_1(D_*)$ is trivial,
we get $\lead_+(\Delta_1(D_*)) = 1$. We get $\lead_+(\Delta_0(D_*)) = 3$ if $\phi(g) \ge 1$,
$\lead_+(\Delta_0(D_*)) = 2$ if $\phi(g) \le -1$, and $5$ if $\phi(g) = 0$. This is a
contradiction.

For another proof based on twisted $L^2$-torsion functions see Subsection~\ref{subsec:The_matrix_3g_plus_2e}.


\typeout{------------------------- Section 10: Examples ------------------------}

\section{Examples}\label{subsection:Examples}



\subsection{Aspherical counterexample and the property special}%
\label{subsec:Aspherical}
  
Hughes--Leary--L\"uck~\cite[Theorem~5.5]{Hughes-Leary-Lueck(2026not_fibre)} constructed
examples of aspherical closed manifolds $M$ in each dimension $d \ge 3$ such that $M$ does
not virtually $\F_2$-fibre, all $L^2$-Betti numbers (over any field) vanish, and the
fundamental group $\pi_1(M)$ is residually (torsionfree and nilpotent) and a
residually  finite Farrell--Jones group. Since $M$ does
not virtually $\F_2$-fibre, there is no finite covering $p \colon N \to M$ for which
there exists a connected finite $CW$-complex $F$ such that $N$ is homotopy equivalent to
the mapping torus $T_f$ for some selfhomotopy equivalence
$f \colon F \xrightarrow{\simeq} F$.

One may ask whether the criterion in
Theorem~\ref{the:Property_special_is_an_obstruction_to_fibring_and_virtually_fibring}
detects this.  The answer is negative in the sense that for any finite covering
$p \colon N \to M$ with connected total space $N$ and for every $\phi \in H^1(N;\IZ)$
there exists a representative $\rho$ of $\overline{\rho}(\widetilde{N};\phi)$ and an integer
$m \in \IZ$ satisfying
\[
  \rho(t) = \begin{cases} 0 & t \le 1;
    \\
   m  \cdot \ln(t) & t \ge 1.
  \end{cases}
\]
Moreover, if $d$ is even, one can additionally arrange $m = 0$.

This can be seen as follows. These manifolds $M$ arise as locally trivial $T^k$-bundles
$M \to B$ such that $B$ is an aspherical closed manifold for $k = 1$ if $\dim(M)$
is odd and $k = 2$ if $\dim(M)$ is even. Let $p \colon N \to M$ be any finite covering and
$\phi$ be any element in $H^1(N;\IZ)$.  Then one can construct a finite covering
$q \colon \widehat{B} \to B$ with a connected total space $\widehat{B} $ and a fibration
$T^k \to N \to \widehat{B}$, see the proof
of~\cite[Theorem~4.11~(iii)]{Hughes-Leary-Lueck(2026not_fibre)}.  Since $\widehat{B}$ is a
closed aspherical manifold, we conclude from Theorem~\ref{the:fibrations_rho}
\[
  \overline{\rho}^{(2)}(\widetilde{N};\phi) = \chi(\widehat{B}) \cdot
  \overline{\rho}^{(2)}(\widetilde{T^k};i^*\phi).
\]
Now $\overline{\rho}^{(2)}(\widetilde{T^1};i^*\phi)$ has a representative $\rho'$ of the
shape
\[
  \rho'(t) = \begin{cases} 0 & t \le 1;
    \\
    m' \cdot \ln(t) & t \ge 1,
  \end{cases}
\]
for some $m' \in \IZ$ by~\cite[(3.23) on page~136]{Lueck(2002)}, whereas
$\overline{\rho}^{(2)}(\widetilde{T^2};i^*\phi)$ vanishes by
Theorem~\ref{the:Basic_properties_of_the_reduced_L2-torsion_function_for_universal_coverings}~%
\ref{the:Basic_properties_of_the_reduced_L2-torsion_function_for_universal_coverings:product_formula}.


\subsection{Twisted torsion of RAAGs}\label{subsec:RAAGs}

Given a finite flag complex $L$, one has the associated RAAG $A_L$. Note that every RAAG
is a residually finite Farrell--Jones group and hence satisfies
Condition~\ref{basic_setup}~\ref{basic_setup:G}.  There is a finite $CW$-complex model for
$BA_L$ coming from the Salvetti complex. If $L$ is $\IQ$-acyclic, $\widetilde{BA_L}$ is
$L^2$-acyclic.  Recall that  the \emph{universal $L^2$-torsion}
\[
  \rho^{(2)}_u(A_L) := \rho^{(2)}_u(\widetilde{BA_L}) \in \Wh^w(A_L)
\]
is defined, see~\cite[Sections~2 and~3]{Friedl-Lueck(2017universal)}.  Let
$\overline{\chi}(\mathrm{Lk}_{L}(v))$ be the reduced Euler characteristic of the link of
the vertex $v \in L^{(0)}$. Let $\mathrm{tor} H_{i-1}(L)$ be the finite subgroup of
$H_{i-1}(L)$ given by the torsion elements. Every vertex $v \in L^{(0)}$ defines a
generator $v \in A_L$, denoted in the same way, and the element $v-1 \in \IZ[A_L]$
determines a class $[v-1]$ in the \emph{weak Whitehead group} $\Wh^w(A_L)$.  For an
integer $m \ge 1$ we denote by $[m] \in \Wh^w(A_L)$ the class of the element
$m \cdot 1_{\IZ[A_L]} \in \IZ[A_L]$.

We will need the following result of
Okun--Schreve~\cite[Corollary~6.4]{Okun-Schreve(2024torsion)}.
(Note that we write the Whitehead group additively in contrast to Okun--Schreve
in~\cite[Corollary~6.4]{Okun-Schreve(2024torsion)}.)

\begin{theorem}\label{thm:uni_torsion_RAAG}
Let $L$ be a $\IQ$-acyclic flag complex.  Then we get in $\Wh^w(A_L)$
\[
\rho^{(2)}_u(A_L)=\sum_{i\geq 2}(-1)^i \cdot [|\mathrm{tor} H_{i-1}(L)|]   +
\sum_{v\in L^{(0)}} \bigl(- \overline{\chi}(\mathrm{Lk}_{L}(v)) \cdot [v-1]\bigr).
\]
\end{theorem}

\begin{theorem}\label{thm:phi_twisted_torsion_RAAGs} 
  Let $L$ be a $\IQ$-acyclic flag complex, let $A_L$ be a right-angled Artin group on $L$, and let $\phi:A_L \to \IZ$ be a group homomorphism.
  Define two constants
\begin{eqnarray*}
  C_L
  & = &
\sum_{i\geq 2}(-1)^i \cdot \ln(|\mathrm{tor} H_{i-1}(L)|);
\\
  C_\phi
  & =&
\sum_{v\in L^{(0)}} |\phi(v)| \cdot \overline \chi(\mathrm{Lk}_{L}(v)).
\end{eqnarray*}
Then a representative $\rho$ for $\overline\rho^{(2)}(A_L;\phi)(t)$ is given by 
\[
t\mapsto
\begin{cases} 
  C_L- C_\phi \cdot \ln(t) & t\geq 1;
  \\
C_L & t \le 1. 
\end{cases}
\]
\end{theorem}
\begin{proof}
  Given a finitely generated group $G$ which is a residually finite Farrell--Jones group,
  and a group homomorphism $\phi \colon G \to \IZ$, there is a well-defined homomorphism
  of abelian groups
  \[
    \zeta_{G;\phi} \colon \Wh^w(G) \to \cale.
  \]
  The homomorphism is defined as follows: we first represent the class $[A] \in \Wh^w(G)$ by a matrix $A \in \Mpar{n}{n}{\IZ[G]}$ for which $r_A^{(2)} \colon L^2(G)^n \to L^2(G)^n$ is a weak isomorphism; then we map $[A]$ to the  twisted
  $L^2$-torsion function $\overline{\rho}^{(2)}(C_*[A];\phi)$ of the $1$-dimensional
  finitely generated free $\IZ[G]$-chain complex $C_*[A]$ whose first differential is
  $r_A \colon \IZ[G]^n \to \IZ[G]^n$. It is uniquely determined by the property that for any
  $L^2$-acyclic finitely generated  free $\IZ[G]$-chain complex $D_*$ the equality
  \[
  \zeta_{G;\phi}(\rho^{(2)}_u(D_*)) = \overline{\rho}^{(2)}(D_*;\phi)
  \]
  holds. Moreover, for an $L^2$-acyclic connected finite $CW$-complex $X$ with $\pi_1(X) = G$
  we obtain 
  \[
  \zeta_{G;\phi}(\rho^{(2)}_u(\widetilde{X})) = \overline{\rho}^{(2)}(\widetilde{X};\phi),
  \]
  see~\cite[Section~3.4]{Friedl-Lueck(2017universal)}.
  
  For $G = \IZ$, the $(1,1)$-matrix $A = (z-1)$ for the standard generator $z \in \IZ$,
  and $m \in \IZ$, the function 
\[
 \rho(t)
  =\begin{cases}
    |m| \cdot \ln(t)  & t\geq 1
    \\
    0 & t\leq 1.
  \end{cases}
\]
is a representative of $\overline{\rho}^{(2)}(C_*[z-1]; m \cdot \id_{\IZ})$. This follows
from combining Theorem~\ref{the:basic_properties_of_overline(rho)(C_ast;varphi)}~%
\ref{the:basic_properties_of_overline(rho)(C_ast;varphi):scaling}
and~\eqref{asymptotics_for_rho_p(t)} with the observations that for $m < 0$ we have
$t^m \ge 1$ if and only if $t \le 1$ and that we can always add $|m| \cdot \ln(t)$ to a
representative.  Thus we get for $g \in A_L$ with $g \not = e$ from
Theorem~\ref{the:basic_properties_of_overline(rho)(C_ast;varphi)}~%
\ref{the:basic_properties_of_overline(rho)(C_ast;varphi):induction} that the function
\[
 \rho(t)
  =\begin{cases}
    |\phi(g)| \cdot \ln(t)  & t\geq 1
    \\
    0 & t \leq 1.
  \end{cases}
\]
is a representative of $\rho^{(2)}(C_*[g-1];\phi) = \zeta_{G;\phi}([g-1])$, for the
$(1,1)$-matrix $(g-1)$ over $\IZ[G]$.  Analogously one easily checks that for an integer
$m \ge 1 $ a representative for $\zeta_{G;\phi}([m])$ is given by the constant function
with value $\ln(m)$ since this function is a representative of
$\overline{\rho}^{(2)}(C_*[m \cdot 1_{\IZ[A_L]}]; \phi)$.  Now
Theorem~\ref{thm:phi_twisted_torsion_RAAGs} follows from
Theorem~\ref{thm:uni_torsion_RAAG}.
\end{proof}

\begin{example}\label{ex:RAAG_RP2}
  Let $L$ be a flag triangulation of the real projective plane on $11$ vertices (see~\cite{Bibby-et-al(2020)} for an example).
  Note   that since $L$ is a triangulation of a manifold, the link of every vertex is homotopic
  to $S^1$.  Whence, for all $v\in L^{(0)}$ we have
  $\overline{\chi}(\mathrm{Lk}_L(v))=-1$.  Let $\phi:A_L\twoheadrightarrow\IZ$ be the
  `Bestvina--Brady character', that is, the map sending each $v$ to $1\in\IZ$.  We have
  that a representative for $\rho^{(2)}(A_L;\phi)$ is given by
  $t\mapsto \ln(2)+11\ln(t)$ for $t\geq 1$ and $t\mapsto\ln(2)$ for $t \leq 1$.  In particular,
  $\rho^{(2)}(A_L;\phi)$ is not $T_0$-special at zero or $T_{\infty} $-special at infinity for any $T_0, T_{\infty} \in\IR_{>0}$.
\end{example}

\begin{example}
  For an abelian group $A$ let $M(A,n)$ be the Moore space of type $(A,n)$.  Let $L$ be a
  flag triangulation of $M(\IZ/2\IZ,2)\vee M(\IZ/2\IZ,3)$.  Now, for any
  $\phi:A_L\twoheadrightarrow \IZ$ there exists an integer constant $C_{\phi,L}$ such that
  $t\mapsto C_{\phi,L}\cdot \ln(t)$ for $t\geq 1$ and $t\mapsto 0$ for $t \leq 1$ is a
  representative for $\rho^{(2)}(A_L;\phi)$.  Indeed, $L$ is $\mathbb{Q}$-acyclic, 
  and the only non-trivial reduced homology groups
    of $L$ are $H_2(L)\cong\mathbb{Z}/2$ and $H_3(L)\cong\mathbb{Z}/2$. Hence
    $C_L = (-1)^3\cdot\ln(2) + (-1)^4\cdot\ln(2) = 0$, and Theorem~\ref{thm:phi_twisted_torsion_RAAGs}
     shows that
    $t \mapsto -C_\phi\cdot\ln(t)$ for $t \ge 1$ and $t\mapsto 0$ for $t\le 1$ is a
    representative of $\overline{\rho}^{(2)}(A_L;\phi)$, where $C_{\phi,L} := -C_\phi$ is
    an integer. In particular $\overline{\rho}^{(2)}(A_L;\phi)$ is special for every
    $\phi$, in contrast to Example~\ref{ex:RAAG_RP2}. 
\end{example}

 
 \subsection{Polynomially growing automorphisms of hyperbolic groups}\label{subsec:poly_growth}
 
 The following theorem is a  reformulation of~\cite[Theorem~3.7]{Hughes-Lueck(2025)}
 taking into account Remark~1.1 of op.\ cit.
 
 \begin{theorem}\label{the:Sum_formula_revisited}
   Let $G$ be a group satisfying Condition~\ref{basic_setup}~\ref{basic_setup:G} and suppose $G$ acts
   cocompactly on a contractible CW-complex $X$ such that the fixed point sets of finite
   subgroups of G are contractible. Suppose that each cell stabilizer $H_\sigma$ of the
   action of G is $L^2$-acyclic and admits a finite model for $BH_\sigma$.
 If  $\phi\colon G \to \IZ$ 
   is a homomorphism, then we get in $\cale$ that
   \[
   \rho^{(2)}(G;\phi) =\sum_{n\geq 0}\sum_{\sigma\in
       I_n}(-1)^n\rho^{(2)}(H_{\sigma};\phi|_{H_\sigma}).
   \]
where $I_n$ is a set of representatives of the $G$-orbits of $n$-cells and $H_{\sigma}$ is the stabilizer of $\sigma$.
 \end{theorem}

 \begin{theorem}\label{the:Gamma_G_rtime_Z}
   Let $\Gamma=G\rtimes_\Theta\IZ$ be a group with $\Gamma$ torsion-free residually finite and $G$
   hyperbolic (or hyperbolic relative to a finite collection of virtually polycyclic
   groups).  Let $\phi\in H^1(\Gamma;\IZ)$.
   If $\Theta$ is polynomially growing, then a
   representative for $\rho^{(2)}(\Gamma;\phi)$ is given by 
   \[
     \rho(t)=\begin{cases}
       0 & \text{if } t\le 1;\\
       k_\phi\ln(t) &\text{if } t\geq 1,
     \end{cases}
   \]
   where $k_\phi$ is a constant depending on $\phi$.
 \end{theorem}
 \begin{proof}
   By~\cite{Andrew-Guerch-Hughes(2026)}, the group $\Gamma$ satisfies the Farrell--Jones
   Conjecture.  In particular, Condition~\ref{basic_setup}~\ref{basic_setup:G} holds for
   $\Gamma$.  The result is proved identically to~\cite[Proposition~7.16]{Hughes-Lueck(2025)} using
   Theorem~\ref{the:Sum_formula_revisited} above.
 \end{proof}

 
 \subsection{Revisiting the matrix 3g+2e}\label{subsec:The_matrix_3g_plus_2e}
 
 The next example is motivated by
 Subsection~\ref{subsec:Back_to_the_non-aspherical_counterexample}. 
 
 Let $G$ be a group satisfying Condition~\ref{basic_setup}~\ref{basic_setup:G}.
 Consider a homomorphism $\phi\colon G\to \IZ$.  Consider $g \in G$ with $|g| = \infty$.
 Let $C_\ast$ be the $1$-dimensional finite free $\IZ[G]$-chain complex whose first
 differential is the map $\IZ[G] \to \IZ[G]$ given by multiplication by $3g+2e$.
 Let $i \colon \IZ \to G$ be the injective homomorphism  sending the
 generator $s \in \IZ$ to $g$. Let $m$ be the integer uniquely determined by
 $\phi \circ i(s) = m$.

We conclude from~\eqref{asymptotics_for_rho_p(t)} and
 Theorem~\ref{the:basic_properties_of_overline(rho)(C_ast;varphi)}~%
\ref{the:basic_properties_of_overline(rho)(C_ast;varphi):induction}
and~\ref{the:basic_properties_of_overline(rho)(C_ast;varphi):scaling}
that a representative for $\rho^{(2)}(C_*;\phi)$ is given
for $m>0$ by
\[
\rho(t) = \begin{cases}
\ln(2) & t \le (2/3)^{1/m};\\
\ln(3) + m\cdot\ln(t) & t \ge (2/3)^{1/m},
\end{cases}
\]
for $m<0$ by
\[
\rho(t) = \begin{cases}
\ln(3) & t \le (3/2)^{1/|m|};\\
\ln(2) + |m|\cdot\ln(t) & t \ge (3/2)^{1/|m|},
\end{cases}
\]
and for $m=0$ by $\rho(t)=\ln(3)$.

This implies that $\rho^{(2)}(C_*;\phi)$ is not $T_0$-special at zero or
$T_{\infty}$-special at infinity for any $T_0, T_{\infty} \in\IR_{>0}$.



\addcontentsline{toc<<}{section}{References} 


\end{document}